\documentclass[11pt]{amsart}
\usepackage{amssymb,amsmath,amsthm}
\usepackage{amsfonts}
\usepackage[dvips]{graphics}
\usepackage{epsfig}
\usepackage{graphicx}
\usepackage{psfrag}
\usepackage{amsmath}
\usepackage{amssymb,subfigure}

\def\R {\mathbb{R}}
\def\N {\mathbb{N}}

\def \and{\quad\text{and}\quad}

\def\ds{\displaystyle}

\def\R{\mathbb R}

\def\N{\mathbb N}
\def\epsilon{\varepsilon}
\def\ds{\displaystyle}

\newcommand{\be}{\begin{equation}}
\newcommand{\ee}{\end{equation}}
\newcommand{\baa}{\begin{array}}
\newcommand{\eaa}{\end{array}}
\newcommand{\ba}{\begin{eqnarray}}
\newcommand{\ea}{\end{eqnarray}}

\numberwithin{equation}{section}

\newtheorem{a3lemma}{Lemma}[section]
\newtheorem{a3theorem}{Theorem}[section]

\newtheorem{a3corollary}{Corollary}[section]
\newtheorem{a3definition}{Definition}[section]

\newtheorem{a3proposition}{Proposition}[section]
\newtheorem{a3remark}{Remark}[section]

\def \no#1#2#3 {{\bf #1} (#3), #2.}
\def \eds#1#2#3 {#1, #2, #3.}
\theoremstyle{definition}

\numberwithin{equation}{section}
\title[On population resilience to external perturbations] {On population resilience to external
perturbations}

\author[L. Roques, and M.D. Chekroun]
{Lionel Roques and Micka\"el D. Chekroun}

\address{Unit\'e Biostatistique et Processus Spatiaux, INRA,
Domaine St Paul - Site Agroparc, 84914 Avignon Cedex 9, France}
\email{lionel.roques@avignon.inra.fr {\rm (L. Roques)}}

\address{\'Ecole Normale Sup\'erieure - CERES-ERTI
\newline\indent
75005 Paris, France} \email{chekro@lmd.ens.fr {\rm (M. D. \
Chekroun)}}

\thanks{The present manuscript has been published as:
\textsc{L. Roques and M.D. Chekroun}, On population resilience to
external perturbations, {\it SIAM Journal on Applied Mathematics
(SIAP)}, {\bf 68} (1), (2007) 133--153. }



\begin{document}
\maketitle

\begin{abstract}
We study a spatially explicit harvesting model in periodic or
bounded environments. The model is governed by a parabolic equation
with a spatially dependent nonlinearity of
Kolmogorov--Petrovsky--Piskunov type, and a negative external
forcing term $-\delta$. Using sub- and supersolution methods and the
characterization of the first eigenvalue of some linear elliptic
operators, we obtain existence and nonexistence results as well as
results on the number of stationary solutions. We also characterize
the asymptotic behavior of the evolution equation as a function of
the forcing term amplitude.

In particular, we define two critical values $\delta^*$ and
$\delta_2$ such that, if $\delta$ is smaller than $\delta^*$, the
population density converges to a ``significant" state, which is
everywhere above a certain small threshold, whereas if $\delta$ is
larger than $\delta_2$, the population density converges to a
``remnant" state, everywhere below this small threshold. Our results
are shown to be useful for studying the relationships between
environmental fragmentation and maximum sustainable yield from
populations. We present numerical results in the case of stochastic
environments.
\end{abstract}

\section{Introduction}\label{section1}

%former Abstract:
%
%We study a spatially explicit harvesting model in periodic or
%bounded environments. The model is governed by a parabolic
%equation with a space-dependent nonlinearity of
%Komolgorov--Petrovsky--Piskunov type, and a negative external
%forcing term. The domain is either the whole space $\mathbb{R}^N$,
%with periodic coefficients, or a bounded domain. Analyzing the
%stationary states, we define two main types of solutions: the
%``significant'' solutions, which always stay above a certain small
%threshold value, and the ``remnant'' solutions, which are always
%below this value. Using sub- and supersolution methods and the
%characterization of the first eigenvalue and first eigenfunction
%of some linear elliptic operators, we obtain existence and
%nonexistence results, as well as results on the number of
%stationary solutions. We also characterize the asymptotic behavior
%of the evolution equation as a function of the forcing term
%amplitude. In particular, we define critical thresholds on the
%forcing term, below which the population density converges to a
%significant state, while it converges to a remnant state whenever
%the forcing term lies above the highest threshold. These bounds
%were shown to be useful in studying the influence of environmental
%fragmentation on the long-time behavior of the population density
%in terms of the forcing term amplitude. We also present numerical
%results in the case of stochastic environments.

Overexploitation has led to the extinction of many
species~\cite{ba}. Traditionally, models of ordinary differential
equations (ODEs) or difference equations have been used to estimate
the maximum sustainable yields from populations and to perform
quantitative analysis of harvesting policies and management
strategies~\cite{getz}. Ignoring age or stage structures as well as
delay mechanisms, which will not be treated by the present paper,
the ODEs models are generally of the type \be \frac{d
U}{dt}=F(U)-Y(U), \label{eqi1} \ee where $U$ is the population
biomass at time $t$, $F(U)$ is the growth function, and $Y(U)$
corresponds to the harvest function. In these models, the most
commonly used growth function is logistic, with $F(U)=U(\mu-\nu U)$
(see \cite{bedmay}, \cite{ms}, \cite{scha}), where $\mu>0$ is the
intrinsic growth rate of the population and $\nu>0$ models its
susceptibility to crowding effects.

Different harvesting strategies $Y(U)$ have been considered in the
literature and are used in practical resource management. A very
common one is the \emph{constant-yield harvesting} strategy, where a
constant number of individuals are removed per unit of time:
$Y(U)=\delta$, with $\delta$ a positive constant. This harvesting
function naturally appears when a quota is set on the
harvesters~\cite{rob2}, \cite{rob1}, \cite{steph}. Another
frequently used harvesting strategy is the \emph{proportional
harvesting} strategy (also called \emph{constant-effort
harvesting}), where a constant proportion of the population is
removed. It leads to a harvesting function of the type $Y(U)=\delta
U$.

Much less has been done in this field using reaction-diffusion
models (but see \cite{shi1}, \cite{neub}, \cite{shi2}). The aim of
this paper is to perform an analysis of some harvesting models,
within the framework of reaction-diffusion equations.

One of the most celebrated reaction-diffusion models was introduced
by Fisher~\cite{fi} and Kolmogorov, Petrovsky, and
Piskunov~\cite{kpp} in 1937 (we call it the Fisher-KPP model). Since
then, it has been widely used to model spatial propagation or
spreading of biological species into homogeneous environments (see
books \cite{ms}, \cite{oku}, and~\cite{tur} for a review). The
corresponding equation is
\begin{equation}\label{eqi5}
u_t=D \nabla^2 u+u(\mu- \nu u),
\end{equation}
where $u=u(t,x)$ is the population density at time $t$ and space
position $x$, $D$ is the diffusion coefficient, and  $\mu$ and $\nu$
still correspond to the \emph{constant} intrinsic growth rate and
susceptibility to crowding effects. In the 1980s, this model was
extended to heterogeneous environments by Shigesada, Kawasaki, and
Teramoto~\cite{skt}. The corresponding model (which we call the
\emph{SKT model} in this paper) is of the type \be\label{eqi4} u_t=
D\nabla ^2 u+u(\mu(x)-\nu(x) u). \ee The coefficients $\mu(x)$ and
$\nu(x)$ now depend on the space variable $x$ and can therefore
include some effects of environmental heterogeneity. More recently,
this model revealed that the heterogeneous character of the
environment plays an essential role in species persistence, in the
sense that for different spatial configurations of the environment a
population can survive or become extinct, depending on the habitat
spatial structure~\cite{bhr1}, \cite{ccL}, \cite{rs1}, \cite{sk}.

As mentioned above, the combination of a harvesting model with a
Fisher-KPP population dynamics model, leading to an equation of the
form $u_t= D \nabla^2 u +u(\mu- \nu u)-Y(x,u)$, has been considered
in recent papers, either using a spatially dependent proportional
harvesting term $Y(x,u)=q(x)u$ in \cite{neub}, \cite{shi2}, or a
spatially dependent and time-constant harvesting term $Y(x)=h(x)$
in~\cite{shi1}. In these papers, the models were considered in
bounded domains with Dirichlet (lethal) boundary conditions.

Here we study a population dynamics model of the SKT type, with a
spatially dependent harvesting term $Y(x,u)$: \be\label{eqi2} u_t= D
\nabla^2 u +u(\mu(x)-\nu(x) u)-Y(x,u). \ee We mainly focus on a
``quasi-constant-yield'' case, where the harvesting  term depends on
$u$ only for very low population densities (ensuring the
nonnegativity of $u$). We consider two types of domains and boundary
conditions. In the first case, the domain is bounded with Neumann
(reflective) boundary conditions; this framework is often more
realistic for modeling species that cannot cross the domain
boundary. In the second case, we consider the model (\ref{eqi2}) in
the whole space $\R^{N}$  with periodic coefficients. This last
situation, though technically more complex, is useful, for instance,
for studying spreading phenomena \cite{bh}, \cite{bhr2}, and for
studying the effects of environmental fragmentation, independently
of the boundary effects. Lastly, note that the effects of
variability in time of the harvesting function will be investigated
in a forthcoming publication~\cite{cr2}.

In section~2, we define a  quasi-constant-yield harvesting
reaction-diffusion model. We prove, on a firm mathematical basis,
existence and nonexistence results for the equilibrium equations, as
well as results on the number of possible stationary states. We also
characterize the asymptotic behavior of the solutions of
(\ref{eqi2}). In section~3, we illustrate the practical usefulness
of the results of section~2, by studying the effects of the
amplitude of the harvesting term on the population density in terms
of environmental fragmentation. Lastly, in section~4, we give new
results for the proportional harvesting case $Y(x,u)=q(x) u$.

\section{Mathematical analysis of a quasi-constant-yield harvesting reaction-diffusion
model}\label{section2}

For the sake of readability, the proofs of the results of section~2
are postponed to section 2.5.

\subsection{Formulation of the model}

In this paper, we consider the model \be\label{eqbevo} u_t= D
\nabla^2 u+u(\mu(x)-\nu(x) u)-\delta h(x)\rho_{\varepsilon}(u),
\quad (t,x) \in \R_+\times\Omega. \ee The function $u=u(t,x)$
denotes the population density at time $t$ and space position $x$.
The coefficient $D$, assumed to be positive, denotes the diffusion
coefficient. The functions $\mu(x)$ and $\nu(x)$ respectively stand
for the spatially dependent intrinsic growth rate of the population,
and for its susceptibility to crowding effects. Two different types
of domains $\Omega$ are considered: either $\Omega=\R^N$ or $\Omega$
is a smooth bounded and connected domain of $\R^N$ ($N \geq 1$). We
qualify the first case as the \textit{periodic case}, and the second
one as the \textit{bounded case}. In the periodic case, we assume
that the functions $\mu(x)$, $\nu(x)$, and $h(x)$ depend on the
space variables in a periodic fashion. For that, let
$L=(L_1,\ldots,L_N)\in (0,+\infty)^{N}$. We recall the following
definition.

\begin{a3definition}
A function $g$ is said to be\/ {\rm L-periodic} if $g(x+k)=g(x)$ for
all $x=(x_1,\ldots,x_N)\in\mathbb{R}^N$ and $k\in
L_1\mathbb{Z}\times\cdots\times L_N\mathbb{Z}$.
\end{a3definition}

Thus, in the periodic case,  we assume that $\mu$, $\nu$, and $h$
are L-periodic. In the bounded case we assume that Neumann boundary
conditions hold: $\frac{\partial u}{\partial n}=0$ on $\partial
\Omega$, where $n$ is the outward unit normal to $\partial \Omega$.
The period cell $C$ is defined by
\[
C:=(0,L_1)\times\cdots\times(0,L_N)
\]
in the periodic case, and in the bounded case we set
\[
C:=\Omega,
\]
for the sake of simplicity of some forthcoming statements.

We furthermore assume that the functions $\mu$ and $\nu$ satisfy
\be\label{hyp} \mu, \nu \in L^{\infty}(\Omega) \quad\text{and}\quad
\exists \ \underline{\nu}\;,\overline{\nu} \in  \R \text{ s.t.~}
0<\underline{\nu}<\nu(x)<\overline{\nu}\quad \forall\  x  \in
\Omega. \ee

Regions with higher values of $\mu(x)$ and lower values of $\nu(x)$
will be qualified as being {\em more favorable}, while, on the other
hand, regions with lower $\mu(x)$ and higher $\nu(x)$ values will be
considered as being {\em less favorable\/} or, equivalently, {\em
more hostile}.

The last term in (\ref{eqbevo}), $\delta h(x)\rho_{\varepsilon}(u)$,
corresponds to a quasi-constant-yield harvesting term. Indeed, the
function $\rho_{\varepsilon}$ satisfies \be\label{hyprho}
\rho_{\varepsilon}\in C^1(\R),\ \rho_{\varepsilon}'\geq0, \
\rho_{\varepsilon}(s)=0 \ \forall s\leq 0 \quad\hbox{and}\quad
\rho_{\varepsilon}(s)=1 \ \forall s\geq \varepsilon, \ee where
$\varepsilon$ is a nonnegative parameter. With such a harvesting
function, the yield is constant in time whenever $u\geq
\varepsilon$, while it depends on the population density when
$u<\varepsilon$. In what follows, the parameter $\varepsilon$ is
taken to be very small. As we prove in the next sections, there are
many situations where the solutions of the model always remain
larger than $\varepsilon$. For these reasons, we qualify our model
as a {\em quasi-constant-yield harvesting SKT model}, the
``dominant'' regime being the constant-yield one. Note that the
function $\rho_\varepsilon$ ensures the nonnegativity of the
solutions of (\ref{eqbevo}). From a biological point of view,
$\varepsilon$ can correspond to a threshold below which harvesting
is progressively abandoned. Considering constant-yield harvesting
functions without this threshold value would be unrealistic since it
would lead to harvest on zero-populations.

Finally, we specify that $\delta \geq 0$ and that $h$ is a function
in $L^{\infty}(\Omega)$ such that \be \exists
  \,\alpha>0 \hbox{ with }\alpha\leq h(x)\leq 1 \,\forall x
\in \Omega. \ee We call $h$ the {\em harvesting scalar field}, and
$\delta$ designates in this way the amplitude of this field.

Before starting our analysis of this model,  we consider the
no-harvesting case, i.e.,~when $\delta=0$. We recall the main known
results in this case. These results will indeed be necessary for the
analysis of the quasi-constant-yield harvesting SKT model.

\subsection{The no-harvesting case}

When $\delta=0$ in (\ref{eqbevo}), our model reduces to the SKT
model described by (\ref{eqi4}). The behavior of the solutions of
this model has been extensively studied in
\cite{bhr1}~and~\cite{bhr2}.

Results are formulated in terms of first (smallest) eigenvalue
$\lambda_1$ of the Schr\"odinger operator $\mathcal{L}_{\mu}$
defined by
\[
\mathcal{L}_{\mu}\phi:=-D \nabla^2-\mu(x)I,
\]
with either periodic boundary conditions (on the period cell $C$) in
the periodic case or Neumann boundary conditions in the bounded
case. This operator is the linearized one of the full model around
the trivial solution. Recall that $\lambda_{1}$ is defined as the
unique real number such that there exists a function $\phi>0$, the
first eigenfunction, which satisfies \be\label{eqblambda1} \left\{
\begin{array}{l@{}}
-D \nabla^2
\phi-\mu(x)\phi=\lambda_{1}\phi\quad \text{in }C,\\
\phi>0\quad\text{in }C,\qquad \|\phi\|_{\infty}=1,
\end{array}
\right. \ee with either periodic or Neumann boundary conditions,
depending on $\Omega$. The function $\phi$ is uniquely defined by
(\ref{eqblambda1})~\cite{bh} and belongs to $W^{2,\tau}(C)$ for all
$1\leq \tau<\infty$ (see \cite{ad} and~\cite{am} for further
details). We set
\[
\underline{\phi}:=\min_{x\in C}\phi(x).
\]

We recall that a stationary state $p$ of (\ref{eqi4}) satisfies the
equation \be\label{eqstabhr1} -D\nabla ^2 p=p(\mu(x)-\nu(x) p). \ee
The following result on the stationary states of (\ref{eqstabhr1})
is proved in~\cite{bhr1}.

\begin{a3theorem} \label{thbhr1}
{\rm(i)} If $\lambda_{1}<0,$ then {\rm(\ref{eqstabhr1})} admits a
unique nonnegative{\em,} nontrivial{\em,} and bounded solution{\em,}
$p_0$.

{\rm(ii)} If $\lambda_{1}\geq 0$, the only nonnegative and bounded
solution of\/ {\rm(\ref{eqstabhr1})} is\/~{\rm0}.
\end{a3theorem}

Moreover, in the periodic case, the solution $p_0$ is L-periodic.
Throughout this paper, $p_0$ always denotes the stationary solution
given by Theorem~\ref{thbhr1}.i.

In order to emphasize that this solution can be ``far'' from~0 (see
Definition~\ref{DefRemn} and the commentary following
(\ref{hypeps0})), we give a lower bound for~$p_0$.

\begin{a3proposition}\label{propbhr1}
Assume that $\lambda_{1}<0;$ then
$p_0\geq\frac{-\lambda_1\underline{\phi}}{\overline{\nu}}$
in\/~$\Omega$.
\end{a3proposition}

The asymptotic behavior of the solutions of (\ref{eqi4}) is also
detailed in~\cite{bhr1}. It is proved that $\lambda_1<0$ is a
necessary and sufficient condition for species persistence, whatever
the initial population $u^0$ is, as follows.

\begin{a3theorem}\label{thdbhr1}
Let $u^{0}$ be an arbitrary bounded and continuous function in\/
$\Omega$ such that $u^0\ge 0,$ $u^0\not\equiv 0$. Let $u(t,x)$ be
the solution of\/ {\rm(\ref{eqi4}),} with initial datum
$u(0,x)=u^0(x)$.

{\rm(i)} If $\lambda_1<0,$ then $u(t,x)\to p_0(x)$ in
$W^{2,\tau}_{loc}\left(\Omega\right)$ for all\/ $1\leq \tau<\infty$
as $t\to +\infty$ (uniformly in the bounded case).

{\rm(ii)} If $\lambda_1\ge 0,$ then $u(t,x)\to 0$ uniformly in\/
$\Omega$ as $t\to+\infty$.
\end{a3theorem}

The situation (i) corresponds to persistence, while in the case~(ii)
the population tends to extinction. In what follows, unless
otherwise specified, we therefore always assume that $\lambda_1<0$,
so that the population survives, at least when there is no
harvesting. We are now in position to start our main analysis of
steady states and related asymptotic behavior of the solutions of
(\ref{eqbevo}).

\subsection{Stationary states analysis}

As is classically demonstrated in finite dimensional dynamical
systems theory and many problems in the infinite dimensional setting
(see, {\it e.g.}, \cite{tem}), the asymptotic behavior of the
solutions of (\ref{eqbevo}) is governed in part by the steady states
and their relative stability properties. In that respect, we study
in this section the positive stationary solutions of (\ref{eqbevo}),
namely the solutions of \be\label{eq_sta_e}
 -D\nabla^2 p_{\delta}=p_{\delta}(\mu(x)-\nu(x)
p_{\delta})-\delta h(x)\rho_{\varepsilon}(p_{\delta}), \quad x \in
\Omega, \ee in the periodic and bounded cases. When needed, we may
write ($\ref{eq_sta_e},\delta$) instead of (\ref{eq_sta_e}).

Note that, provided $p_{\delta}\geq\varepsilon$ in $\Omega$,
$p_\delta$ is equivalently a solution of the simpler equation
\be\label{SIAPeq_sta} -D\nabla^2 p_{\delta}=p_{\delta}(\mu(x)-\nu(x)
p_{\delta})-\delta h(x), \quad x \in \Omega. \ee This last equation
has been analyzed in the case of Dirichlet boundary conditions
in~\cite{shi2}, in the particular case of constant coefficients
$\mu$ and~$\nu$.

Because of the type of harvesting function considered here, we are
led to introduce the following definition.

\begin{a3definition}\label{DefRemn}
Set $\varepsilon_0:= \frac{\varepsilon}{\underline{\phi}}\geq
\varepsilon$. We say that a nonnegative function $\sigma$ is\/ {\rm
remnant} whenever $\max_{C} \sigma<\varepsilon_0,$ whereas it is
{\rm significant} if it is a bounded function satisfying
$\min_{C}\sigma \geq \varepsilon_0$.
\end{a3definition}

\begin{figure}[t!]
\centering
\includegraphics[width=8cm,height=5cm]{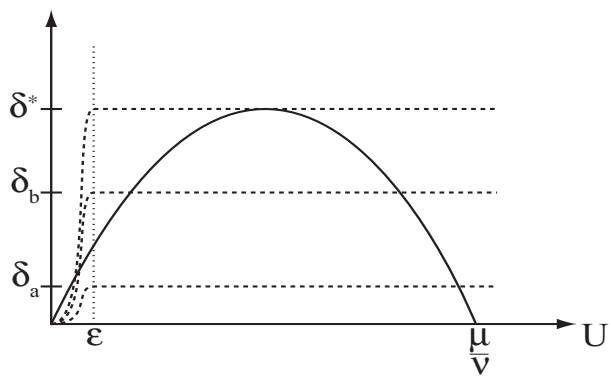}
\caption{The logistic growth function $U\mapsto U(\mu-\nu U)$ (solid
line){\em,} and the harvesting function $U\mapsto \delta
\rho_\varepsilon(U)$ for three values of $\delta$ (dashed lines).
The abscissae of the points of intersection of the solid and dashed
lines correspond{\em,} respectively{\em,} to remnant (if smaller
than $\varepsilon$) and significant (if strictly larger than
$\varepsilon$) stationary solutions of\/ {\rm(\ref{ode})}. We
observe that the number of significant solutions is as follows: one
if $\delta<k(\varepsilon)$ (case $\delta=\delta_a$){\em;} two if
$k(\varepsilon)\leq \delta<\mu^2/(4 \nu)$ (case
$\delta=\delta_b$){\em;} one if $\delta= \mu^2/(4 \nu)$ (case
$\delta=\delta^*$){\em;} zero if $\delta>\mu^2/(4 \nu) $. The number
of nonzero remnant solutions is zero or more if $\delta \leq
k(\varepsilon)$ (depending on the shape of $\rho_\varepsilon$){\em;}
one or more if $\delta > k(\varepsilon),$ since, from\/
{\rm(\ref{hyprho}),} $\rho^{\prime}_\varepsilon(0)=0$. We assumed
here that $\varepsilon_0=\varepsilon$.}\label{fig:ode}
\end{figure}

%{\em Remark 1:}
\begin{a3remark}\label{SIAPrem1}\rm
The concepts of  remnant and significant solutions, as well as the
harvesting term $\delta h(x) \rho_\varepsilon(u)$, are not
classical. In order to clarify these notions, we present in
Figure~\ref{fig:ode} a short graphical study of the nonspatial model
\be \frac{dU}{dt}=U(\mu-\nu U)-\delta \rho_\varepsilon(U)=:k(U),
\quad t\in\R_+, \label{ode} \ee with constant coefficients
$\mu,\nu>0$.
\end{a3remark}

Since $\varepsilon_0$ is assumed to be small in our model, the
remnant solutions of (\ref{eq_sta_e}) correspond to very low
population densities. On the other hand, significant solutions are
everywhere above $\varepsilon_0$. In particular, a constant yield is
ensured in that case. In contrast to the ODE case, stationary
solutions which are neither remnant nor significant may exist, as
outlined in the next theorems. However, as we will see while
studying the long-time behavior of the solutions of the model
(\ref{eqbevo}), they are of less importance (see Theorem~\ref{th2}
and section~3). The threshold $\varepsilon_0$ is different from
$\varepsilon$ in general. We had to define remnant and significant
functions using $\varepsilon_0$ for technical reasons (see the proof
of Theorem~\ref{th1}.ii, equation (\ref{useeps0})). Since
$\varepsilon$ is assumed to be very small, it has no implication on
the biological interpretation of our results. Moreover, most of our
results still work when $\varepsilon_0$ is replaced
by~$\varepsilon$.

Let us now start our analysis of (\ref{eq_sta_e}). In what follows,
we always assume that \be \varepsilon_0 <
\ds{\frac{-\lambda_1\underline{\phi}}{4
\overline{\nu}}},\label{hypeps0} \ee so that, in particular, from
Proposition~\ref{propbhr1}, the solution $p_0$ of (\ref{eqstabhr1})
is significant.

We begin by proving that there exists a threshold $\delta^{\ast}$
such that, if the amplitude $\delta$ is below $\delta^{\ast}$,
(\ref{eq_sta_e}) admits significant solutions, while it does not in
the other case.

\begin{a3theorem}\label{th0}
Assume that $\lambda_1<0;$ then there exists $\delta^{\ast}\geq0$
such that

{\rm(i)} if $\delta\leq \delta^{\ast},$ there exists at least a
positive significant solution $p_\delta\leq p_0$ of\/
{\rm(\ref{eq_sta_e});}

{\rm(ii)} if $\delta>\delta^{\ast},$ there  is no positive
 significant solution of\/ {\rm(\ref{eq_sta_e})}.
\end{a3theorem}

%{\em Remark 2:}
\begin{a3remark}\label{rem2}\rm
There  is no positive bounded solution of (\ref{eq_sta_e}) whenever
$\lambda_1\geq 0$.
\end{a3remark}

Under stronger hypotheses, we are able to prove that
(\ref{eq_sta_e}) admits \textit{at most} two significant solutions.
In order to state this result, we need some definitions. Let $G$ be
the space  defined by \be\label{defG2} G:=H^1(C) \ee in the bounded
case, and by
 \be \label{defG1}
 G:=H^1_{per}=\left\{\psi
\ \in \ H^{1}_{loc}(\mathbb{R}^N) \text{ such that } \psi \text{ is
L-periodic}\right\}
 \ee
in the periodic case. Let us define the standard Rayleigh quotient:
for all $\psi\in G$, $\psi\not\equiv 0$, and for all $\sigma \in
L^{\infty}(C)$, \be
\mathcal{R_\sigma(\psi)}:=\frac{{\int_{C}}D|\nabla\psi|^2-\sigma(x)
\psi^{2}}{{\int_{C}}\psi^{2}}. \label{defR} \ee According to the
Courant--Fischer theorem (see, e.g.,~\cite{bbbs}), the second
smallest eigenvalue $\lambda_2$ of the operator $\mathcal{L}_\mu$
can be characterized by \be \lambda_2=\min_{E_k\subset G,
\text{dim}(E_k)=2}\max_{\psi \in E_k, \ \psi\not\equiv 0}
\mathcal{R_\mu}(\psi). \label{minmax} \ee This characterization is
equivalent to the classical one given in~\cite{gt}.

We are now in position to state the following theorem.

\begin{a3theorem}\label{th1.1}
Assume that $\lambda_1<0\leq\lambda_2;$ then{\em,} in the bounded
case{\em,} {\rm(\ref{eq_sta_e})} admits at most two significant
solutions. In the periodic case{\em,} {\rm(\ref{eq_sta_e})} admits
at most two L\/$^{\prime}$-periodic significant solutions for all
L\/$^\prime\in(0,+\infty)^N$. Moreover{\em,} under these
hypotheses{\em,} if two solutions $p_{1,\delta}$ and $p_{2,\delta}$
exist{\em,} they are ordered in the sense that{\em,} for
instance{\em,} $p_{1,\delta}<p_{2,\delta}$ in\/~$\Omega$.
\end{a3theorem}

%{\em Remark 3:}
\begin{a3remark}\label{rem3}\rm
Similar methods also allow us to assess a result on the number of
solutions of (\ref{SIAPeq_sta}). Indeed, if
$\lambda_1<0\leq\lambda_2$, then we obtain that (\ref{SIAPeq_sta})
admits at most two nonnegative bounded (and periodic in the periodic
case) solutions. If these solutions exist, they are ordered.
\end{a3remark}

In the periodic case, Theorem~\ref{th1.1} also gives some
information on the periodicity of the significant solutions of
(\ref{eq_sta_e}), which are actually found to have the same
periodicity as the coefficients of (\ref{eq_sta_e}), as seen in the
next result.

\begin{a3corollary}\label{cor1.1}
Assume that $\lambda_1<0\leq\lambda_2$.  Then, in the periodic case,
the significant periodic solutions of\/  {\rm(\ref{eq_sta_e})} are
L-periodic.
\end{a3corollary}

The fact that  $\lambda_1<0$ is directly related to the instability
of the trivial solution in the SKT model. The additional condition
$\lambda_2\geq 0$ in this theorem is linked to the existence of a
stable manifold or center manifold of the steady state~0 of the SKT
model, in some appropriate functional spaces (see~\cite{tem}).
Therefore, the assumptions of Theorem~\ref{th1.1}, and the Krein
Rutmann theory, allow us to conclude that under these assumptions
the unstable manifold of~0 is of dimension equal to \textit{one} or
equivalently the stable manifold is of codimension~1. Such results
on multiplicity of solutions of elliptic nonlinear equations with a
source or sink term have been investigated in the past and are known
nowadays as being of Ambrosetti-problem type. These results also
involve manifolds of codimension~1 (in the functional space of
forcing) and first and second eigenvalues (for the Laplace operator
only) (see~\cite{nir} for a survey of these results).

In any event, Theorem~\ref{th1.1} relies on the assumption that
$\lambda_2 \geq 0$. In the next proposition, we give conditions
under which $\lambda_2$ may become positive.

\begin{a3proposition}
{\rm(i)} In the bounded case{\em,} if $C$ is a (smooth) domain with
diameter $d:=\max_{x,y\in C}\|x-y\|_{\R^N},$ $\lambda_2(C)\geq
D(\frac{\pi}{d})^2-\max_C \mu$.

{\rm(ii)}~In the periodic case{\em,} $\lambda_2(C)\geq
D(\frac{\pi}{L_d})^2-\max_C \mu,$ where $L_d$ denotes the length of
the longest diagonal of the period cell~$C$.\label{proplbl2}
\end{a3proposition}

For instance, when $C=[0,1]\times[0,1]$, we have $d=L_d=\sqrt{2}$;
thus, for $D=1$ and $\max_C{\mu}=4$, we get $\lambda_2>0.9$.
However, this lower bound is far from being optimal. Indeed, in all
our computations of section~\ref{sec_appl}, and under the same
hypothesis on $C$ and $D$, we always had $\lambda_2>0$, while
$\max_C{\mu}=10$. Sharper lower bounds for $\lambda_{2}$ can be
found in~\cite{bp}; however, those bounds are also more sensitive to
the geometry of the domain and thus less general. They are therefore
not detailed here.

We now introduce a result which is important for more applied
ecological questions. Indeed, one of the main drawbacks of
Theorem~\ref{th0} is that it gives no computable bound for
$\delta^{\ast}$. Obtaining information on the value of
$\delta^{\ast}$ is precious for  ecological questions such as the
study of the relationships between $\delta^{\ast}$ and the
environmental heterogeneities. The next theorem states some
computable estimates~of~$\delta^{\ast}$.

Let us define \be \delta_1:=\ds{\frac{\lambda_1^2
\underline{\phi}}{\overline{\nu}(1+\underline{\phi})^2}}
\quad\hbox{and}\quad \delta_2:=\ds{\frac{\lambda_1^2 }{4 \alpha
\underline{\nu}}}. \label{defd1d2} \ee Note that neither $\delta_1$
nor $\delta_2$ depend on $\delta$~and~$\varepsilon$.

\begin{a3theorem}\label{th1}
{\rm(i)} If $\lambda_1<0$ and $\delta\leq \delta_1,$ then there
exists a positive significant (and $L$-periodic in the periodic
case) solution $p_\delta$ of\/ {\rm(\ref{eq_sta_e})} such that
$p_\delta \geq -\frac{\lambda_1
\phi}{\overline{\nu}(1+\underline{\phi})}$.

{\rm(ii)} If $\lambda_1<0$ and $\delta> \delta_2,$ the only possible
positive bounded  solutions of\/ {\rm(\ref{eq_sta_e})} are remnant.
\end{a3theorem}

The lower bound of part~(i), for $p_\delta$, does not depend on
$\varepsilon$. Thus, there is a clear distinction between the
remnant and significant solutions. Note that, of course,
$\delta_1\leq\delta_2$.

The formulae (\ref{defd1d2}) allow numerical evaluations. An
important quantity to compute is the size of the gap
$\delta_2-\delta_1$ and its fluctuations
 in  terms of environmental configurations. This question is
addressed in section~\ref{sec_appl} through a numerical study.

\subsection{Asymptotic behavior}\label{sub_asymp}

In this section, we prove that the quantity $\delta^{\ast}$ in fact
corresponds to a maximum sustainable yield, in the sense that when
$\delta$ is smaller than $\delta^{\ast}$, the population density
$u(t,x)$ converges to a significant stationary state of
(\ref{eqbevo}) as $t\to \infty$, whereas when $\delta$ is larger
than $\delta^{\ast}$, the population density converges to a
stationary state which is not significant. In fact, when $\delta$ is
larger than the quantity $\delta_2$ defined by (\ref{defd1d2}) we
even prove that the population converges to a remnant stationary
state of (\ref{eqbevo}).

We assume here that the harvesting starts on a stabilized population
governed by the standard SKT model with $\delta=0$. From
Theorem~\ref{thdbhr1}, this means that we study the behavior of the
solutions $u(t,x)$ of our model (\ref{eqbevo}), starting with the
initial datum $u(0,x)=p_0(x)$. Since we have assumed that
$\lambda_1<0$, it follows from Theorem~\ref{thbhr1},
Proposition~\ref{propbhr1}, and (\ref{hypeps0}) that $p_0$ is well
defined and significant.

Let us describe, with the next theorem, the long-time behavior of
the population density.

\begin{a3theorem}
Let $u(t,x)$ be the solution of {\rm(\ref{eqbevo})} with initial
datum $u(0,x) = p_0(x)$. Then $u$ is nonincreasing in $t$ and the
following hold:

{\rm(i)} If $\delta\leq \delta^{\ast},$ $u(t,x)\to p_\delta(x)$
uniformly in\/ $\Omega$ as $t\to +\infty,$ where $p_{\delta}$ is the
maximal significant  solution of\/ {\rm(\ref{eq_sta_e})}.
Moreover{\em,} $p_\delta$ is L-periodic in the periodic case.

{\rm(ii)} If $\delta > \delta^{\ast},$ then the function
$u(t,\cdot)$ converges uniformly in\/ $\Omega$ to a solution of\/
{\rm(\ref{eq_sta_e})} which is not significant.

{\rm(iii)} If $\delta > \delta_2,$  the function $u(t,\cdot)$
converges uniformly in\/ $\Omega$ to a remnant solution of\/
{\rm(\ref{eq_sta_e})}.\label{th2}
\end{a3theorem}

%{\em Remark 4:}
\begin{a3remark}\label{rem4}\rm
If, in addition, we assume that $\lambda_2\geq 0$, then
Theorem~\ref{th1.1} says that, whenever $\delta\leq \delta^{\ast}$,
(\ref{eqbevo}) admits at most two significant stationary states
(which are periodic stationary states in the periodic case). In that
case, the stationary state $p_\delta$ selected at large times is the
higher one. If we do not assume that $\lambda_2\geq 0$, this
stationary state can still be defined as ``the maximal one'' that
can be constructed by a sub- and supersolution method
(see~\cite{am2}).
\end{a3remark}

From the above theorem, we observe that, whenever
$\delta\leq\delta^{\ast}$, the solution $u(t,x)$ of (\ref{eqbevo}),
with initial datum $p_0$, remains significant for all times $t \geq
0$. This ensures a constant yield in time and justifies the name of
the model.

Similar results could be obtained for a wider class of initial data.
Indeed, with similar methods, the convergence of $u(t,x)$ to a
significant solution of (\ref{eq_sta_e}) can be obtained whenever
$\delta\leq\delta^{\ast}$ for all bounded and continuous initial
data $u(0,x)$ which are larger than the smallest significant
solution of (\ref{eq_sta_e}). In particular, when $u(0,x)$ is larger
than the maximal significant solution of (\ref{eq_sta_e}), $u(t,x)$
converges to this maximal significant solution as $t\to +\infty$. A
more detailed analysis of the basin of attraction related to the
maximal significant solution will be further investigated in the
forthcoming paper~\cite{cr2}.

Theorem~\ref{th2} shows that the practical determination of
$\delta^*$ is directly linked to the size of the gap
$\delta_{2}-\delta_{1}$. As we will see in section~\ref{sec_appl},
this gap $(\delta_1,\delta_2)$ can be very narrow in certain
situations. In those cases, the numerical computation of $\delta_1$
and $\delta_2$ therefore gives a sharp localization of the maximum
sustainable quota $\delta^{\ast}\in[\delta_1,\delta_2]$, which can
be of nonnegligible ecological interest.

\subsection{Proofs of the results of section~2}\mbox{}

{\it Proof of Proposition\/ {\rm\ref{propbhr1}}}. Let $\phi$ be
defined by (\ref{eqblambda1}), with the appropriate boundary
conditions. Set $\kappa_0:=\frac{-\lambda_1}{\overline{\nu}}$. Then
the function $\kappa_0 \phi$ satisfies
\begin{eqnarray*}
-D\nabla^2(\kappa_0\phi)-\mu(x) \kappa_0\phi +\nu(x)
(\kappa_0\phi)^2 &
= & \lambda_1 \kappa_0 \phi +\nu(x) (\kappa_0\phi)^2 \\
&= & \kappa_0 \phi (\lambda_1 +\nu(x) \kappa_0\phi)\leq 0.
\end{eqnarray*}
Thus $\kappa_0\phi$ is a subsolution of (\ref{eqstabhr1}) satisfied
by $p_0$. Since for $M\in\R$ large enough $M$ is a supersolution of
(\ref{eqstabhr1}), it follows from the uniqueness of the positive
bounded solution $p_0$ of (\ref{eqstabhr1}) that $p_0\geq
\kappa_0\phi \geq\frac{-\lambda_1\underline{\phi}}{\overline{\nu}}$.
\qquad $\Box$

Before proving Theorem~\ref{th0}, we begin with the following lemma.

\begin{a3lemma}
For all $\delta>0,$ if $p_\delta$ is a nonnegative bounded solution
of {\rm(\ref{eq_sta_e}),} then $p_\delta\leq p_0$.\label{lem1}
\end{a3lemma}

{\it Proof of Lemma\/ {\rm\ref{lem1}}}.  Assume that there exists
$x_0\in\Omega$ such that $p_\delta(x_0)> p_0(x_0)$. The function
$p_\delta$ satisfies
\[
-D\nabla^2p_\delta-p_\delta(\mu(x)-\nu(x)p_\delta)=-\delta h(x)
\rho_\varepsilon(p_\delta)\leq 0,
\]
and thus $p_\delta$ is a subsolution of (\ref{eqstabhr1}) satisfied
by $p_0$. Since for $M\in\R$ large enough $M$ is a supersolution of
(\ref{eqstabhr1}), we can apply a classic iterative method to infer
the existence of a solution $p^{\prime}_{0}$ of (\ref{eqstabhr1})
(with Neumann boundary conditions in the bounded case since both
$p_\delta$ and $M$ satisfy Neumann boundary conditions) such that
$p_\delta\leq p^{\prime}_0\leq M$. In particular,
$p^{\prime}_0(x_0)>p_0(x_0)$, which is in contradiction with the
uniqueness of the positive bounded solution of (\ref{eqstabhr1}).
\qquad $\Box$
%\endproof %

{\it Proof of Theorem\/ {\rm\ref{th0}}}. Let us define
\[
\delta^{\ast}:=\sup\{\delta\geq 0, {\rm(\ref{eq_sta_e})} \hbox{
admits a significant solution}\}.
\]
For $\delta=0$, we know from Proposition~\ref{propbhr1} that $p_0$
is a significant solution of (\ref{eq_sta_e}). Moreover, for
$\delta$ large enough, the nonexistence of significant solutions of
(\ref{eq_sta_e}) is a direct consequence of the maximum principle
(it is also a consequence of the proof of Theorem~\ref{th1}.ii).
Thus $\delta^{\ast}$ is well defined and bounded.

Assume that $\delta^{\ast}>0$, and let us prove that
($\ref{eq_sta_e},\delta^{\ast}$) admits a significant solution. By
definition of $\delta^{\ast}$, there exists a sequence
$(p_{\delta_k})_{k\in\N}$ of solutions of
($\ref{eq_sta_e},\delta_k$) with  $0<\delta_k\leq\delta^{\ast}$ and
$\delta_k\to\delta^{\ast}$ as $k\to +\infty$. Moreover, from
Lemma~\ref{lem1}, $\varepsilon_0\leq p_{\delta_k}\leq p_0$ for all
$k\geq0$. Thus, from standard elliptic estimates and Sobolev
injections, the sequence $(p_{\delta_k})_{k\in\N}$ converges (up to
the extraction of some subsequence) in $W^{2,\tau}_{loc}$, for all
$1\leq\tau<\infty$, to a significant solution $p_{\delta^{\ast}}$ of
($\ref{eq_sta_e},\delta^{\ast}$).

Now, let $0\leq\delta<\delta^{\ast}$. Then
\[
-D\nabla^2p_{\delta^{\ast}}-p_{\delta^{\ast}}(\mu(x)-\nu(x)p_{\delta^{\ast}})+\delta
h(x)=(\delta-\delta^{\ast})h(x)<0,
\]
and thus $p_{\delta^{\ast}}$ is a subsolution  of
($\ref{eq_sta_e},\delta$). Since $p_0$ is a supersolution of
($\ref{eq_sta_e},\delta$), and $p_{\delta^*}\leq p_0$, a classical
iterative method gives the existence of a significant solution
$p_\delta$ of ($\ref{eq_sta_e},\delta$) (with Neumann boundary
conditions in the bounded case since both $p_0$ and $p_\delta$
satisfy Neumann boundary conditions). This concludes the proof of
Theorem\linebreak \ref{th0}. \qquad $\Box$
%\endproof

{\it Proof of Theorem\/ {\rm\ref{th1.1}}}.  As a preliminary, we
prove that if two solutions exist, then they cannot intersect. Let
$p_{1,\delta}$ and $p_{2,\delta}$ be two significant solutions of
(\ref{eq_sta_e}). In the bounded case, we assume that $p_{1,\delta}$
and $p_{2,\delta}$ satisfy Neumann boundary conditions.  In the
periodic case, we assume that there exists L$^{\prime}\in
(0,+\infty)^N$ such that $p_{1,\delta}$ and $p_{2,\delta}$ are
L$^\prime$-periodic, and then denote the period cell by
$C^{\prime}$. Let us set $q_\delta:=p_{2,\delta}-p_{1,\delta}$. Then
$q_\delta$ verifies \be -D \nabla^2 q_\delta
-[\mu(x)-\nu(x)(p_{1,\delta}+p_{2,\delta})]q_\delta=0; \ee thus,
setting $\rho(x):=\mu(x)-\nu(x)(p_{1,\delta}+p_{2,\delta})$, we
obtain \be \label{eqqd} -D \nabla^2 q_\delta -\rho(x) q_\delta=0,
\ee with the same boundary conditions that were satisfied by
$p_{1,\delta}$~and~$p_{2,\delta}$.

Let $\widehat{\lambda_1}$ and $\widehat{\lambda_2}$ be respectively
the first and second eigenvalues of the operator
$\mathcal{L}_\rho:=-D \nabla^2-\rho I$. Let
$\mathcal{R}_\sigma(\phi)$, be defined by (\ref{defR}). Since
$\rho(x)<\mu(x)-2\underline{\nu}\varepsilon_0$ for all $x\in
\Omega$, we get
\[
\mathcal{R_\rho(\varphi)}\geq \mathcal{R_\mu(\varphi)}+2
\underline{\nu}\varepsilon_0
\]
for all $\varphi\in G^{\prime}$, where $G^{\prime}:=H^1(C)$ in the
bounded case and
\[
G^{\prime}:=H^1_{per}=\left\{\varphi \ \in \
H^{1}_{loc}(\mathbb{R}^N) \text{ such that } \varphi \text{ is
 L}^{\prime}\hbox{-periodic}\right\}
 \]
in the periodic case. Thus, by the classical min-max formula
(\ref{minmax}), it follows that \be \widehat{\lambda_2}\geq
\lambda_2 +2 \underline{\nu}\varepsilon_0>0. \label{seceig} \ee
Furthermore, from (\ref{eqqd}), 0~is an eigenvalue of the operator
$\mathcal{L}_\rho$. Thus, (\ref{seceig}) implies that
$\widehat{\lambda_1}=0$. As a consequence, $q_\delta$ is a principal
eigenfunction of the operator $\mathcal{L}_\rho$. The principal
eigenfunction characterization thus implies that $q_\delta$ has a
constant sign. Finally, we get that $p_{1,\delta}$ and
$p_{2,\delta}$ do not intersect each other.

Let us now prove that (\ref{eq_sta_e}) admits at most two
significant solutions. Arguing by contradiction, we assume that
there exist three significant (L$^{\prime}$-periodic in the periodic
case, for some L$^{\prime}\in (0,+\infty)^N$) solutions
$p_{1,\delta}$, $p_{2,\delta}$, and $p_{3,\delta}$ of
(\ref{eq_sta_e}). From the above result, we may assume, without loss
of generality, that
$p_{3,\delta}>p_{2,\delta}>p_{1,\delta}>\varepsilon_0$. Set
$q_{2,1}:=p_{2,\delta}-p_{1,\delta}$  and
$q_{3,2}:=p_{3,\delta}-p_{2,\delta}$; then these functions satisfy
the equations \be \label{eqqd2} -D\nabla^2 q_{2,1} -\rho_{2,1}(x)
q_{2,1}=0 \ee and \be \label{eqqd3} -D\nabla^2 q_{3,2}
-\rho_{3,2}(x) q_{3,2}=0, \ee with
$\rho_{2,1}:=\mu(x)-\nu(x)(p_{1,\delta}+p_{2,\delta})$ and
$\rho_{3,2}:=\mu(x)-\nu(x)(p_{2,\delta}+p_{3,\delta})$.  Moreover,
$q_{2,1}>0$ and $q_{3,2}>0$. Thus 0~is the first eigenvalue of the
operators $\mathcal{L}_{\rho_{2,1}}:=-D \nabla^2-\rho_{2,1}I$ and
$\mathcal{L}_{\rho_{3,2}}:=-D \nabla^2-\rho_{3,2}I$ with either
Neumann or L$^{\prime}$-periodic boundary conditions.

From the strong maximum principle (see, e.g.,~\cite{gt}) (together
with Hopf's lemma in the bounded case, and using the
L$^{\prime}$-periodicity of $q_{3,2}$ in the periodic case), we
obtain the existence of $\theta>0$ such that $q_{3,2}>\theta$. Since
the operator $\mathcal{L}_{\rho_{3,2}}$ is self-adjoint, we have the
following formula for its first eigenvalue
$\widehat{\lambda_1}^{3,2}$:
\[
\widehat{\lambda_1}^{3,2}=\min_{\varphi \in
G^{\prime}}\mathcal{R}_{\rho_{3,2}}(\varphi).
\]
Thus
\begin{eqnarray*}
\widehat{\lambda_1}^{3,2}=\min_{\varphi \in
G^{\prime}}\left\{\mathcal{R}_{\rho_{2,1}}(\varphi)+\frac{\int_C\nu(p_{3,\delta}-p_{1,\delta})\varphi^2}{\int_C\varphi^2}\right\}&
\geq & \min_{\varphi \in
G^{\prime}}\left\{\mathcal{R}_{\rho_{2,1}}(\varphi)\right\}+\underline{\nu}\theta\\
&\geq & \widehat{\lambda_1}^{2,1} + \underline{\nu}\theta,
\end{eqnarray*}
where $\widehat{\lambda_1}^{2,1}$ is the first eigenvalue of the
operator $\mathcal{L}_{\rho_{2,1}}$. Since the first eigenvalues of
the operators $\mathcal{L}_{\rho_{2,1}}$ and
$\mathcal{L}_{\rho_{3,2}}$ are both~0, we deduce that $0\geq 0+
\underline{\nu}\theta >0$, hence a contradiction. \qquad $\Box$

{\it Proof of Corollary\/ {\rm\ref{cor1.1}}}. Let $p_\delta$ be a
significant L$^\prime$-periodic solution of (\ref{eq_sta_e}), and
let $k\in \prod_{i=1}^{N}L_{i}\mathbb{Z}$. From the L-periodicity of
(\ref{eq_sta_e}), $p_\delta({\cdot} +k)$ is also a solution of
(\ref{eq_sta_e}). By periodicity  of $p_\delta$, the functions
$p_\delta$ and $p_\delta(\cdot+k)$ intersect each other. Thus, from
Theorem~\ref{th1.1}, since  $p_\delta$ and $p_\delta(\cdot+k)$ are
both L$^\prime$-periodic, $p_\delta\equiv p_\delta(\cdot+k)$.
Therefore, $p_\delta$ is an L-periodic function. \qquad $\Box$

{\it Proof of Proposition\/ {\rm\ref{proplbl2}}}.  In the bounded
case, let $\tilde{C}$ be the convex hull of the set $C$. It was
proved in~\cite{pw} that the second Neumann eigenvalue of the
Laplace operator $-D \nabla^2$ on $\tilde{C}$ was larger than
$D(\frac{\pi}{d})^2$. Since $C\subset\tilde{C}$, we have
$H^1(C)\subset H^1(\tilde{C})$. Using formula (\ref{minmax}), we
thus obtain that the second eigenvalue of $\mathcal{L}_\mu$ in the
bounded case satisfies $\lambda_2\geq D(\frac{\pi}{d})^2 -\max_C
\mu$. This proves part~(i) of Proposition~\ref{proplbl2}.

In the periodic case, since $H^1_{per}$ can be seen as a subset of
$H^1(C)$, it follows from (\ref{minmax}) that \be \lambda_2\geq
\min_{E_k\subset H^1(C), \text{dim}(E_k)=2}\max_{\psi \in E_k, \
\psi\not\equiv 0} \mathcal{R_\mu}(\psi). \label{mima2} \ee The
period cell $C$ is convex but not smooth enough to assert that the
right-hand side of (\ref{mima2}) is equal to the second eigenvalue
in the bounded case. Let $L_d$ be the longest diagonal of $C$. Then
$C$ is included in a ball $B_{L_d}$ of diameter $L_d$. Thus, from
formula (\ref{minmax}), the right-hand side of (\ref{mima2}) is
larger than the second eigenvalue of $\mathcal{L}_\mu$ on $B_{L_d}$.
From~(i), the conclusion of~(ii) follows. \qquad $\Box$

{\it Proof of Theorem\/~{\rm\ref{th1}}, part\/ {\rm(i)}}. Let
$\lambda_1$ and $\phi$ be defined by (\ref{eqblambda1}), and let
$\kappa$ be a nonnegative real number such that $\kappa >
\varepsilon_0$.  Then we have \be \baa{r@{\,}c@{\,}l@{}} -D\nabla^2
(\kappa \phi)-\kappa \phi (\mu(x)-\kappa
\phi \nu(x))+\delta h(x)\rho_{\varepsilon}(\kappa \phi) & \leq & \lambda_1 \kappa \phi + \kappa^2 \phi^2 \nu(x) +\delta  \\
& \leq &\kappa \phi (\lambda_1+\kappa\phi \nu(x))+\delta \\
& \leq & \ds{\max_{\tau\in I}\{\tau(\lambda_1+\tau
\overline{\nu})\}+\delta }, \eaa \label{eqbp1.1} \ee where $I= \{
\kappa \phi(x)$, $x  \in  C \}$. Setting
$g(\tau):=\tau(\lambda_1+\tau \overline{\nu})$, since
$\|\phi\|_{\infty}=1$, and since $g$ is a convex function, it
follows from (\ref{eqbp1.1}) that \be -D\nabla^2(\kappa \phi)-\kappa
\phi (\mu(x)-\kappa \phi \nu(x))+\delta
h(x)\rho_{\varepsilon}(\kappa \phi) \leq
\max\{g(\kappa),g(\kappa\underline{\phi})\}+\delta . \label{eqbp1.2}
\ee Let us take $\kappa_0$ be such that
$g(\kappa_0)=g(\kappa_0\underline{\phi})$, namely $\kappa_0=
-\frac{\lambda_1}{\overline{\nu}(1+\underline{\phi})}$ (note that
$\kappa_0 \phi > \varepsilon$). We get \be -D\nabla^2(\kappa_0
\phi)-\kappa_0 \phi (\mu(x)-\kappa_0 \phi \nu(x))+\delta h(x)\leq
\ds{-\frac{\lambda_1^2
\underline{\phi}}{\underline{\nu}(1+\underline{\phi})^2}}+\delta
\leq 0, \label{eqbp1.3} \ee from the hypothesis on $\delta$ of
Theorem~\ref{th1}.i. Therefore, $\kappa_0 \phi$ is a subsolution of
(\ref{eq_sta_e}) with either L-periodic or Neumann boundary
conditions. Moreover, if $M$ is a large enough constant, $M$ is a
supersolution of (\ref{eq_sta_e}) with L-periodic or Neumann
boundary conditions. Thus, it follows from a classical iterative
method that there exists a solution $\underline{p_\delta}$ of
(\ref{eq_sta_e}), with the required boundary conditions, and which
satisfies $\kappa_0\phi \leq \underline{p_\delta} \leq M$ in
$\Omega$. Moreover, in the periodic case, since $\kappa_0 \phi$ and
$M$ are L-periodic and since (\ref{eq_sta_e}) is also L-periodic, it
follows that $p_\delta$ is L-periodic. Theorem~\ref{th1}.i
is\linebreak proved. \qquad $\Box$

{\it Proof of Theorem\/ {\rm\ref{th1}}, part\/ {\rm(ii)}}. Assume
that $\lambda_1<0$, $\delta>\delta_2$, and that there exists a
positive bounded solution $p_\delta$ of (\ref{eq_sta_e}) which is
not remnant; i.e., \be \exists\ x_0 \ \text{with} \
p_\delta(x_0)\geq\varepsilon_0. \label{eqbnonex1} \ee Since $\phi$
is bounded from below away from~0 and $p_\delta$ is bounded, we can
define
\begin{equation}
\gamma^{*}=\inf\left\{\gamma>0, \ \gamma\phi>p_\delta \text{ in
}\Omega\right\}> 0. \label{eqbp1.5}
\end{equation}
It follows from the definition of $\gamma^*$ that $\gamma^*\phi\geq
p_\delta$ in $\Omega$, and in particular, $\gamma^*\phi(x_0)\geq
p_\delta(x_0)\geq\varepsilon_0$. Since
 $\|\phi\|_\infty=1$, we get
$\gamma^*\geq \varepsilon_0$. Thus, \be \gamma^*\phi\geq
\varepsilon_0\underline{\phi}=\varepsilon, \label{useeps0} \ee which
implies $\rho_{\varepsilon}( \gamma^* \phi)=1$. Thus,
$h(x)\rho_{\varepsilon}( \gamma^* \phi)\geq \alpha$, and we get
\[
-D\nabla^2(\gamma^* \phi)-\gamma^* \phi (\mu(x)-\gamma^* \phi
\nu(x))+\delta h(x) \rho_{\varepsilon}( \gamma^* \phi)  \geq
\gamma^* \phi (\lambda_1+\gamma^*\phi \nu(x))+\delta \alpha
\]
on $\Omega$. Moreover, since $\gamma^*\phi>0$ and $\nu\geq
\underline{\nu}$, we have $\gamma^* \phi (\lambda_1+\gamma^*\phi
\nu(x))\geq  -\smash{\frac{\lambda_1^2}{4\underline{\nu}}}$. Using
the fact that $\delta>\delta_2= \frac{\lambda_1^2 }{4 \alpha
\underline{\nu}}$, we thus get \be -D\nabla^2(\gamma^*
\phi)-\gamma^* \phi (\mu(x)-\gamma^* \phi \nu(x))+\delta h(x)
\rho_{\varepsilon}( \gamma^* \phi)  \geq
\ds{-\frac{\lambda_1^2}{4\underline{\nu}}+\delta\alpha}>0
\label{sur1} \ee on $\Omega$. Therefore, $\gamma^* \phi$ is a
supersolution of (\ref{eq_sta_e}). Set $z:=\gamma^{*}\phi-p_\delta$.
From the definition of $\gamma^*$, we know that $z\geq 0$ and that
there exists a sequence $(x_{n})_{n\in\N}$ in $\Omega$ such that
$z(x_{n})\to 0$ as $n \to +\infty$.

In the bounded case, up to the extraction of some subsequence,
$x_{n}\to \overline{x}  \in  \Omega$ as $n\to +\infty$. By
continuity, $z(\overline{x})=0$. Moreover, subtracting
(\ref{eq_sta_e}) from (\ref{sur1}), we get \be\label{eqz1}
-D\nabla^2 z +[\nu(x)(\gamma^*\phi+p_{\delta})+\chi(x)-\mu(x)]z >0
\quad\text{in }\Omega, \ee where the function $\chi$ is defined by
$\chi(x)=\delta h(x) \frac{\rho_{\varepsilon}( \gamma^*
\phi(x))-\rho_{\varepsilon}( p_{\delta}(x))}{\gamma^*
\phi(x)-p_{\delta}(x)}$ whenever $\gamma^* \phi(x)-p_{\delta}(x)\neq
0$, and $\chi(x)=\rho_{\varepsilon}^{\prime}(p_{\delta}(x))$
otherwise. Since $\rho_\varepsilon$ is $C^1$, $\chi$ is bounded.
Thus $b(x):=\nu(x)(\gamma^*\phi+p_{\delta})+\chi(x)-\mu(x)$ is a
bounded function. Using the strong elliptic maximum principle, we
deduce from (\ref{eqz1}) that $z\equiv 0$. Thus $\gamma^{*}\phi
\equiv p_\delta$ is a positive solution of (\ref{eq_sta_e}). It is
in contradiction with (\ref{sur1}).

In the periodic case, we must also consider the situation where the
sequence $(x_n)_{n\in\N}$ is not bounded. Let $(\overline{x}_{n})
\in \overline{C}$ be such that $x_{n}-\overline{x}_{n} \in
\prod_{i=1}^{N}L_{i}\mathbb{Z}$. Up to the extraction of some
subsequence, we can assume that there exists $\overline{x}_{\infty}
\in \overline{C}$ such that $\overline{x}_{n}\to
\overline{x}_{\infty}$ as $n\to +\infty$. Set
$\phi_n(x)=\phi(x+x_n)$ and $p_{\delta,n}(x)=p_\delta(x+x_{n})$.
From standard elliptic estimates and Sobolev injections, it follows
that (up to the extraction of some subsequence) $p_{\delta,n}$
converge in $W^{2,\tau}_{loc}$, for all $1\leq \tau<\infty$, to a
function $p_{\delta,\infty}$ satisfying
\[
-\nabla^2(D
p_{\delta,\infty})-p_{\delta,\infty}(\mu(x+\overline{x}_{\infty})-p_{\delta,\infty}\nu(x+\overline{x}_{\infty}))+\delta
h(x+\overline{x}_{\infty})\rho_{\varepsilon}(p_{\delta,\infty})=0
\]
in $\R^N$, while $\gamma^*\phi_n$ converges to
$\gamma^*\phi_{\infty}:=\gamma^*\phi(\cdot+\overline{x}_{\infty})$,
and
\[
-\nabla^2(D\gamma^*\phi_\infty)
-\gamma^*\phi_\infty(\mu(x+\overline{x}_{\infty})-\gamma^*\phi_\infty\nu(x+\overline{x}_{\infty}))+\delta
h(x+\overline{x}_{\infty})\rho_{\varepsilon}(\gamma^* \phi_\infty)>0
\]
in $\R^N$. Let us set
$z_{\infty}(x):=\gamma^*\phi_{\infty}(x)-p_{\delta,\infty}(x)$. Then
$z_{\infty}(x)=\lim_{n\to +\infty}z(x+x_{n})$, and therefore
$z_{\infty}\geq 0$ and $z_{\infty}(0)=0$. Moreover, there exists a
bounded function $b_\infty$ such that \be\label{eqz2} -D\nabla^2
z_\infty +b_\infty z_\infty >0 \quad\hbox{in }\R^N. \ee It then
follows from the strong maximum principle that $z_{\infty}\equiv 0$,
and we again obtain a contradiction. Finally, we necessarily have
$p_\delta\leq \varepsilon_0$, and the proof of Theorem~\ref{th1}.ii
is complete. \qquad $\Box$

{\it Proof of Theorem\/ {\rm\ref{th2}}, part\/ {\rm(i)}}.
 Assume that $\delta\leq
\delta^{\ast}$. Let $p_\delta$ be the unique maximal significant
solution defined in the proof of Theorem~\ref{th1}.i. Then, from
Lemma~\ref{lem1}, \be p_\delta (x) \leq p_0(x) = u(0,x) \quad
\forall x \in \Omega, \label{eqdp_pd} \ee which implies \be
p_\delta(x)\leq u(t,x) \quad\text{in }\mathbb{R}_{+}\times\Omega,
\label{eqdv<u} \ee since $p_\delta$ is a stationary solution of
(\ref{eqbevo}). Moreover, since $p_0$ is a supersolution of
(\ref{eq_sta_e}), $u$ is nonincreasing in time $t$, and standard
parabolic estimates imply that $u$ converges in
$W_{loc}^{2,\tau}\left(\Omega\right)$, for all $1\leq\tau<\infty$,
to a bounded stationary solution $u_{\infty}$ of (\ref{eqbevo}).
Furthermore, from (\ref{eqdv<u}) we deduce that
 $p_\delta \leq u_{\infty}\leq p_0$. Since $p_\delta$ is
 the maximal positive solution of (\ref{eq_sta_e}), it follows that $u_\infty\equiv
 p_\delta$. Moreover, in the periodic case, since $p_0$ and
 (\ref{eqbevo}) are L-periodic, $u(t,x)$ is also L-periodic in $x$. Therefore the convergence is uniform in
 $\Omega$. Part~(i) of Theorem~\ref{th2} is proved.
\qquad $\Box$

{\it Proof of Theorem\/~{\rm\ref{th2}}, parts\/ {\rm(ii)}
and\/~{\rm(iii)}}. Assume that $\delta > \delta^{\ast}$. Since 0~is
a stationary solution of (\ref{eqbevo}) and $u(0,x)=p_0>0$, we
obtain that $u(t,x)>0$ in $\R^+\times\Omega$, and again, from
standard parabolic estimates, we know that $u$ converges in
$W_{loc}^{2,\tau}\left(\Omega\right)$ (for all $1\leq \tau<\infty$)
to a bounded stationary solution $\underline{u_{\infty}}\geq 0$ of
(\ref{eqbevo}) as $t\to +\infty$. Moreover, in the periodic case,
from the L-periodicity of the initial data and of (\ref{eqbevo}), we
know that $u(t,\cdot)$ and $\underline{u_{\infty}}$ are L-periodic.
Therefore the convergence is uniform in $\Omega$. It follows from
Theorem~\ref{th0}.ii that $\underline{u_{\infty}}$ cannot be a
significant solution of (\ref{eq_sta_e}). Moreover, if
$\delta>\delta_2$, Theorem~\ref{th1}.ii ensures that
$\underline{u_{\infty}}$ is a remnant solution of (\ref{eq_sta_e}).~
\qquad $\Box$

\section{Numerical investigation of the effects of environmental fragmentation}\label{sec_appl}

We propose here to apply the results of section 2, on the estimation
of the maximum sustainable yield, to the study of the effects of
environmental fragmentation.  A theoretical investigation of the
relationships between maximum sustainable yield and fragmentation is
difficult to achieve (see Remark~\ref{rem5}). To overcome this
difficulty, we propose a numerical study in the case of stochastic
environments. First, we show that the gap $\delta_2-\delta_1$,
obtained from (\ref{defd1d2}) and Theorem~\ref{th1}, remains small
whatever the degree of fragmentation is. This gap corresponds to the
numerical values of the harvesting quota $\delta$ for which we do
not know whether the population density will converge to a
significant or a remnant solution of the stationary equation
(\ref{eq_sta_e}). Second, we show that there is a monotone
increasing relationship between the maximal sustainable yield
$\delta^{\ast}$ and the habitat aggregation.

%{\em Remark 5:}
\begin{a3remark}\label{rem5}\rm
In a periodic environment, a simple way  of changing the degree of
fragmentation without changing the relative spatial pattern
(favorable area/unfavorable area ratio) is to modify the size of the
period cell $C$. Assume that $\mu(x)=\eta(\frac{x}{L})$, for some
$1$-periodic function $\eta$ with positive integral and for some
$L>0$. This means that the environment consists of square cells of
side $L$. Setting $\lambda_{1,L}:=\lambda_1$ and $\phi_L:=\phi$, we
then have $-D\Delta \phi_L -\eta\left(\frac{x}{L}\right)
\phi_L=\lambda_{1,L}\phi_L$ on $[0,L]^N$. The function
$\psi_L(x):=\phi_L(L x)$ thus satisfies $-D\Delta
\psi_L-L^2\eta(x)\psi_L=L^2 \lambda_{1,L}\psi_L$ in $[0,1]^N$, with
1-periodicity. From the Rayleigh formula we thus obtain
\[
\lambda_{1,L}=\min_{\psi \in H^1_{per}}
\frac{D}{L^2}\frac{\int_{[0,1]^N}|\nabla\psi|^2}{\int_{[0,1]^N}\psi^2}-\frac{\int_{[0,1]^N}\eta
\psi^2}{\int_{[0,1]^N}\psi^2};
\]
therefore $\lambda_{1,L}<0$ (since $\psi \equiv 1 \in H^1_{per}$),
and $\lambda_{1,L}$ decreases with $L$. It implies that $\delta_2$
increases with $L$. The relationship between $\delta_1$ and $L$ is
less clear since $\ds \underline{\phi_L}=\min\nolimits_{C} \phi_L$
may not always be an increasing function of~$L$.
\end{a3remark}

In order to lessen the boundary effects and to focus on
fragmentation, we place ourselves in the periodic case. For our
numerical computations, we assume that the environment is made of
two components, favorable and unfavorable regions. This is expressed
in the model (\ref{eqbevo}) through the coefficient $\mu(x)$, which
takes two values $\mu^+$ or $\mu^-$, depending on the space variable
$x$. We also assume that
\[
\mu^+>\mu^-, \quad\nu(x)\equiv 1, \quad h(x)\equiv 1, \hbox{ and
}D=1.
\]

\begin{figure}[t!]
\centering
\includegraphics[width=8cm]{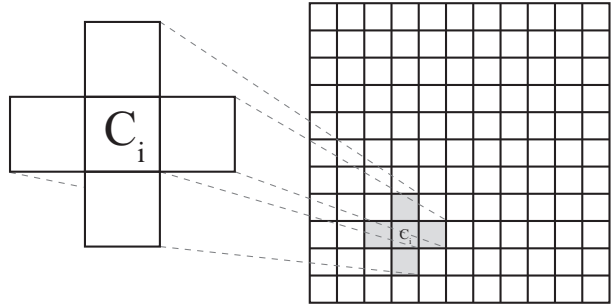}
\caption{The\/ {\rm4}-neighborhood system{\em:} an element $C_i$ of
$C$ and its four neighbors.}\label{fig:vois4}
\end{figure}

\begin{figure}[t!]
\centering
\subfigure[$s=3400$]{\label{figp+e}\includegraphics[width=3.5cm]{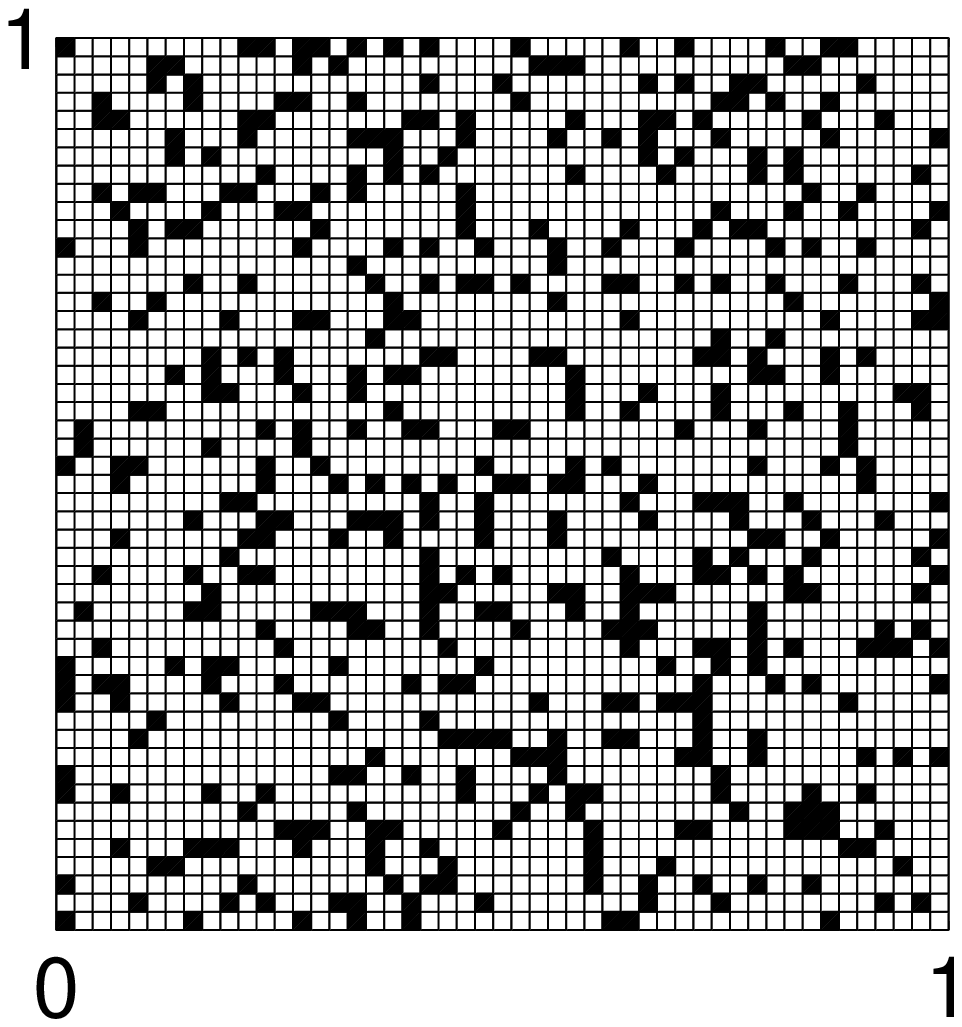}}
\subfigure[$s=3800$]{\label{figp+g}\includegraphics[width=3.5cm]{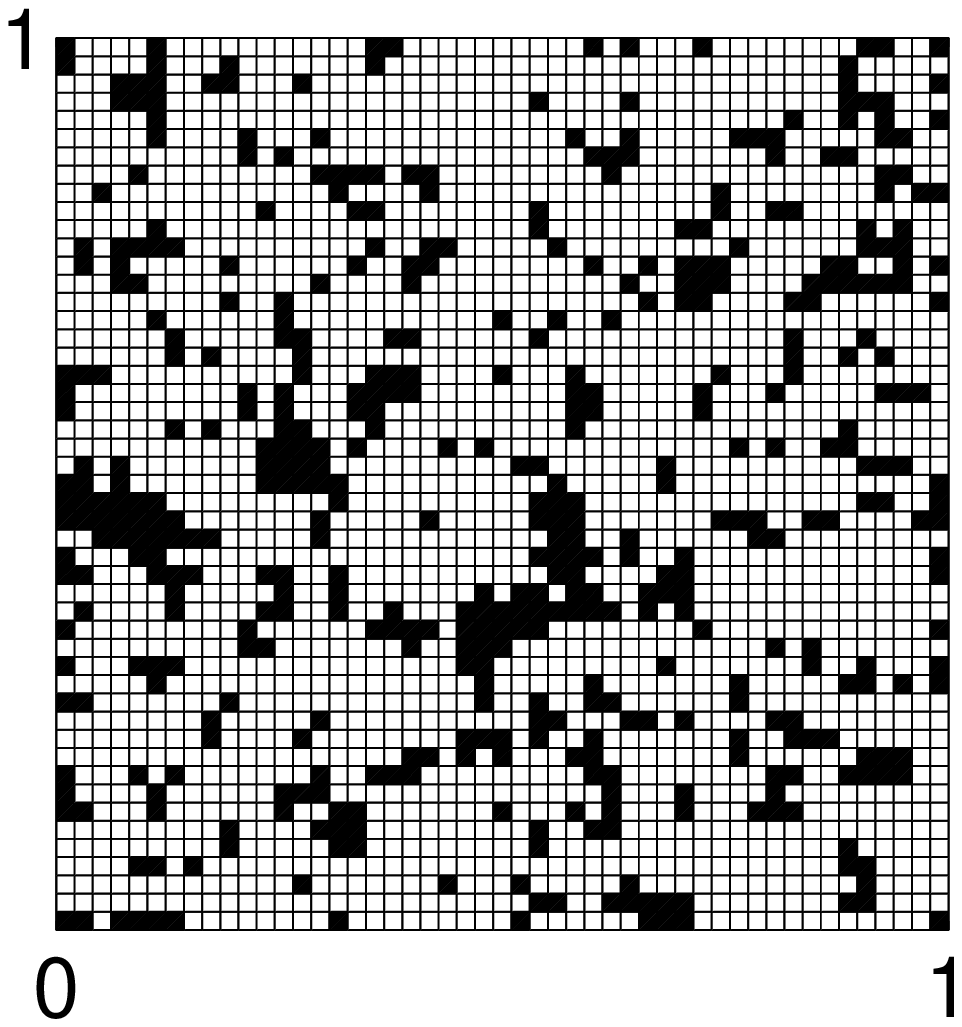}}
\subfigure[$s=4200$]{\label{figp+i}\includegraphics[width=3.5cm]{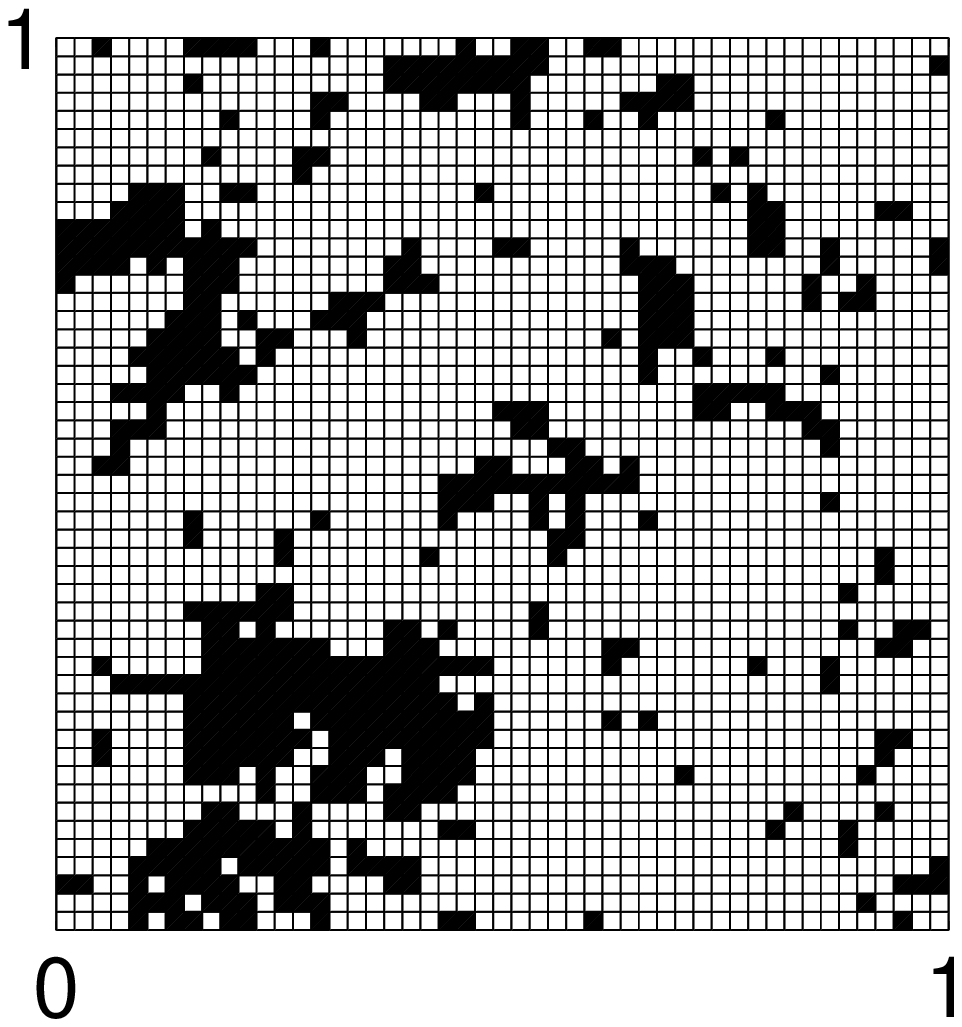}}
\subfigure[$s=4600$]{\label{figp+j}\includegraphics[width=3.5cm]{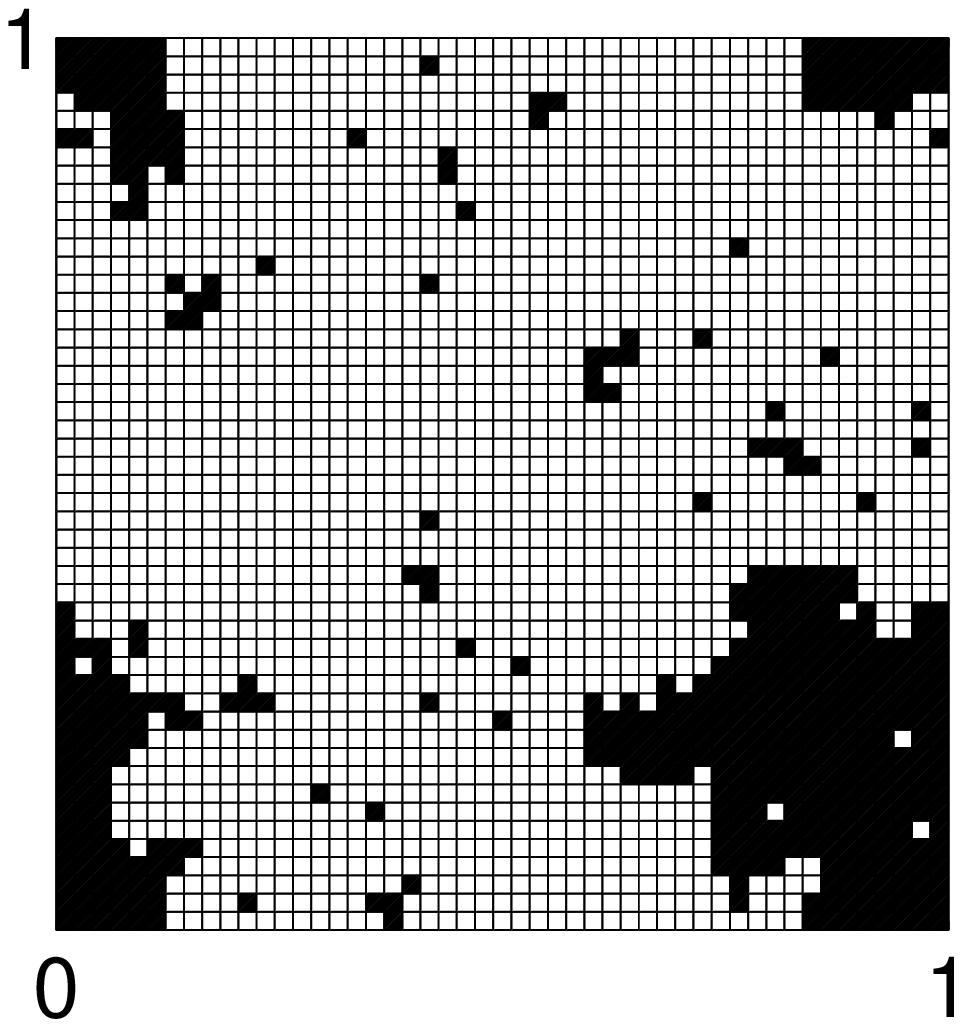}}
\subfigure[$s=4800$]{\label{figp+nj}\includegraphics[width=3.5cm]{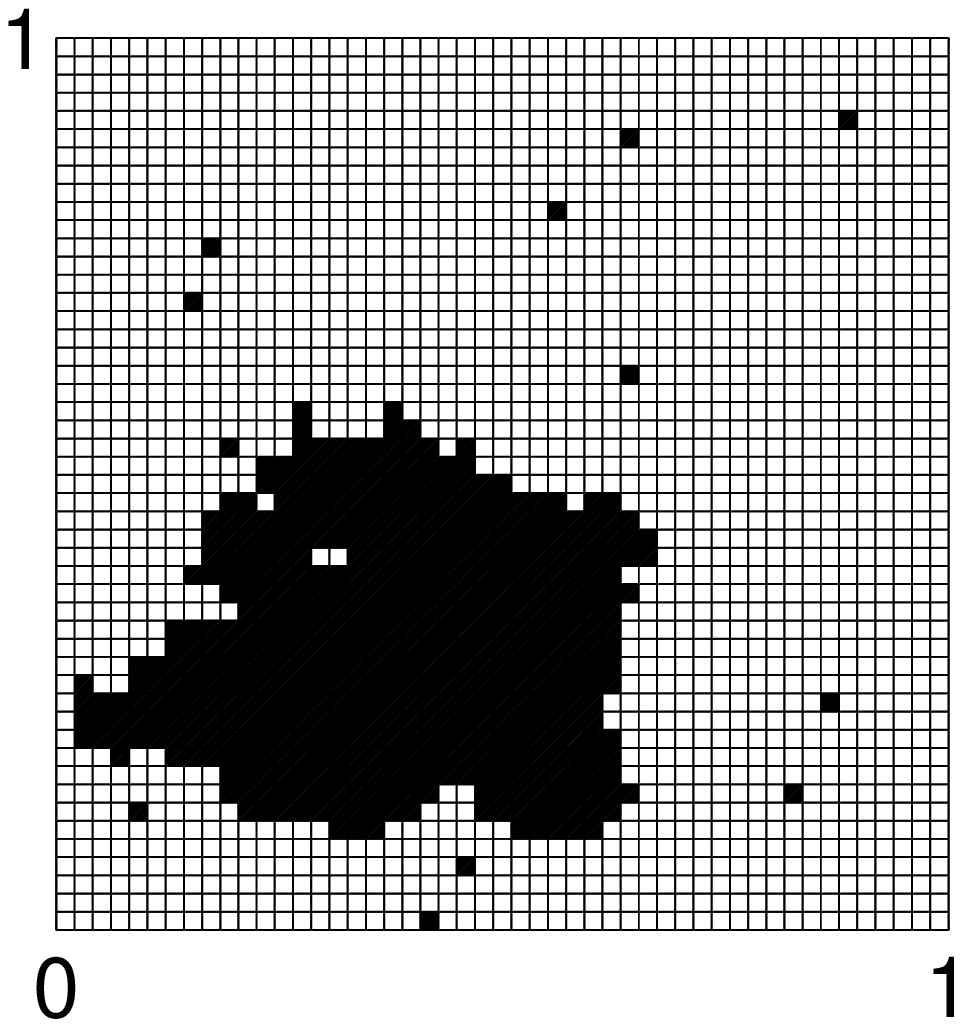}}
\subfigure[$s=4900$]{\label{figp+k}\includegraphics[width=3.5cm]{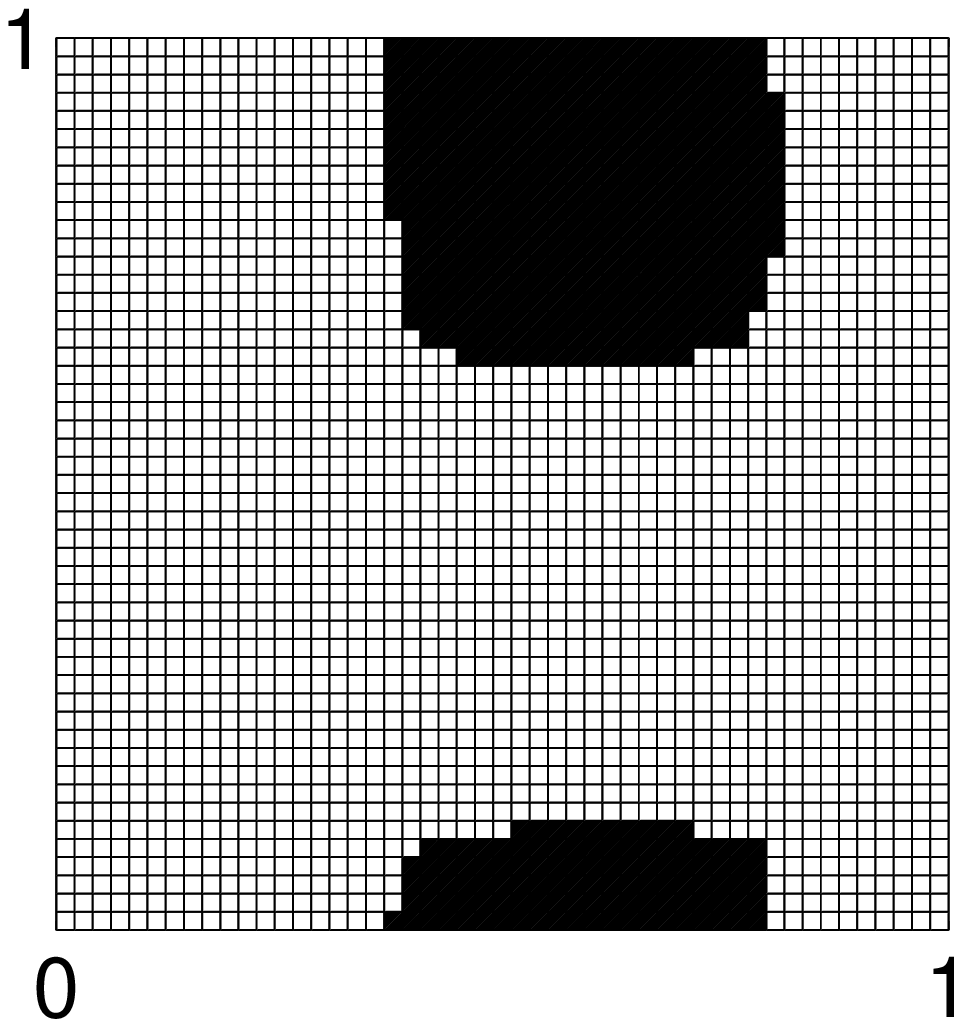}}
\caption{Some samples of the landscapes used for the computations of
$\delta_1$ and $\delta_2,$ with different values of the habitat
aggregation index $s$. The black areas correspond to more favorable
environment{\em,} where $\mu(x)=\mu^+$.}\label{SIAPfig:frag}
\end{figure}

Using a stochastic model for landscape generation \cite{rs1}, we
built 2000 samples of binary environments, on the two-dimensional
period cell  $C=[0,1]^2$, with different degrees of fragmentation.
In all
 these environments, the favorable region, where $\mu(x)=\mu^+$,
 occupies 20\%~of
 the period cell. The environmental
 fragmentation is defined as follows. We discretize the cell $C$
 into $n_C=50\times 50$ equal squares $C_i$. The lattice made of
the cells $C_i$ is equipped with a 4-neighborhood system $V(C_i)$
(see Figure~\ref{fig:vois4}), with toric conditions. On each cell
$C_i$, we assume that the function $\mu$ takes either the value
$\mu^+$ or $\mu^-$, while the number $n_{+}=\hbox{card}\{i$,
$\mu\equiv \mu^+$ on $C_i\}$ is fixed to $n_C \times
\frac{20}{100}=500$.
 For each landscape sample $\omega=(\mu(C_i))_{i=1,\ldots, n_C}$, we set
$s(\omega)=\frac{1}{2}\sum_{C_i\subset C}\sum_{C_j \in V(C_i)}
1\!\!1\{\mu(C_j)=\mu(C_i)\}$,  the number of pairs of
 neighbors $(C_i,C_j)$ such that $\mu$ takes the same value on $C_i$ and $C_j$ ($1\!\!1\{\cdot\}$
 is the indicator function). The number $s(\omega)$ is
 directly linked to the environmental fragmentation: a landscape pattern
 is all the more aggregated as $s(\omega)$ is high, and all the
 more fragmented as $s(\omega)$ is small
(Figure \ref{SIAPfig:frag}). Thus, we shall refer to $s$ as the
``habitat aggregation index.''

%{\em Remark 6:}
\begin{a3remark}\label{rem6}\rm
There exist several ways of obtaining hypothetical landscape
distributions. The commonest  are neutral landscape models,
originally introduced by Gardner {\it et al.} \cite{gard}. They can
include parameters which regulate the fragmentation \cite{keitt}. We
preferred to use a stochastic landscape model presented
in~\cite{rs1}, since it allows an exact control of the favorable and
unfavorable surfaces and is therefore well adapted for analyzing the
effects of fragmentation per se. This model is inspired from
statistical physics. The number of pairs of similar neighbors $s$ is
controlled during the process of landscape generation. This quantity
can be measured a posteriori on the landscape samples. Other
measures of fragmentation could have been used, such as fractal
dimension (see~\cite{mandel}). For a discussion on the different
ways of measuring habitat fragmentation in real-world situations,
the interested reader can refer to \cite{SIAPfahrig2003}.
\end{a3remark}

\begin{figure}[t!]
\centering
\includegraphics[width=10cm]{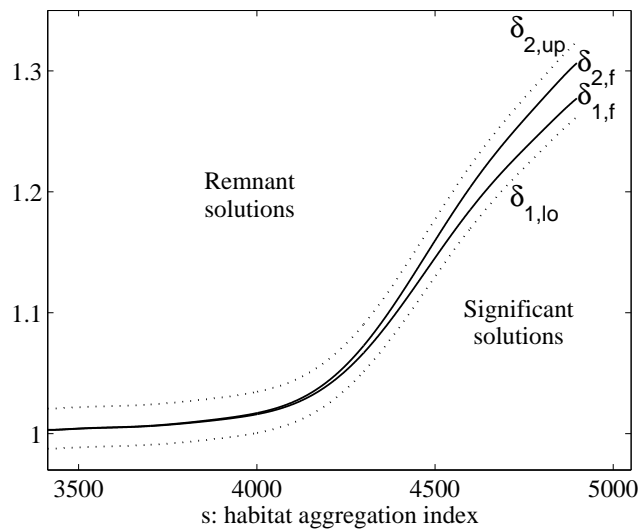}
\caption{Solid lines{\em:} $\delta_{1,f}$ and $\delta_{2,f}$
correspond respectively to the data sets\/
$\{(s^{i},\delta_1^{i})\}_{i=1,\ldots, 2000}$  and\/
$\{(s^{i},\delta_2^{i})\}_{i=1,\ldots, 2000},$ fitted with ninth
degree polynomials. Dashed lines{\em:} $\delta_{1,lo}$ is a lower
prediction bound for new observations of $\delta_{1},$ and
$\delta_{2,up}$ an upper prediction bound for new observations of
$\delta_{2},$ with in both cases a certainty level
of\/~{\rm99\%}.}\label{fig:d1etd2}\vspace*{-1pc}
\end{figure}

For our computations, we took $\mu^+=10$ and $\mu^-=0$, and we
computed the corresponding values of $\lambda_{1}^{i}$,
$\delta_1^{i}$, and $\delta_2^{i}$ on each landscape sample
$\omega^{i}$ of aggregation index $s^{i}$, for $i=1,\ldots, 2000$.
The eigenvalues $\lambda_1^{i}$ were computed with a finite elements
method. We fitted the data sets
$\{(s^{i},\delta_1^{i})\}_{i=1,\ldots, 2000}$  and
$\{(s^{i},\delta_2^{i})\}_{i=1,\ldots, 2000}$ using ninth degree
polynomials (it is enough to assess whether the relations between
$s$ and $\delta_{1}, \delta_{2}$ tend to be monotonic or not). The
resulting fitted curves $\delta_{1,f}$ and $\delta_{2,f}$ are
presented in Figure~\ref{fig:d1etd2}. Under the assumption of
normally distributed values of $\delta_1$ and $\delta_2$ for fixed
$s$ values, we computed a lower prediction bound ($\delta_{1,lo}$)
for new observation of $\delta_1$ and an upper prediction bound for
$\delta_2$ ($\delta_{2,up}$), with a level of certainty of~99\%.
Thus, given a configuration $\omega$, with a fixed value of $s$,
when $\delta$ is smaller than $\delta_{1,lo}$ we take a 0.5\%~chance
of being above $\delta_1$, while when $\delta$ is larger than
$\delta_{2,up}$ we take a  0.5\%~chance of being below $\delta_2$.
The small thickness of the intervals $(\delta_{1,lo},\delta_{2,up})$
emphasizes the quality of the relationship between the habitat
aggregation index $s$ and the maximum sustainable yield
$\delta^{\ast}\in [\delta_1,\delta_2]$. This also indicates that the
criteria of Theorems \ref{th1} and~\ref{th2} are close to being
optimal, at least in some situations.

Furthermore, as we can observe, the values of $\delta_1$ and
$\delta_2$ tend to increase as $s$ increases, and thus as the
environment aggregates. Since $\delta^{\ast}\in
[\delta_1,\delta_2]$, we deduce from the computations presented in
Figure~\ref{fig:d1etd2} that $\delta^{\ast}$ tends to increase with
environmental aggregation.

These tests were performed for particular values of $\mu^+$ and
$\mu^-$. However,  the thickness of the interval
$(\delta_1,\delta_2)$ can be determined for all values of $\mu^{+},
\mu^-$ without further numerical computations, provided that
$\mu^+-\mu^-=10$. Indeed, let us set $B:=\mu^+-\mu^-$. For a fixed
value of $B$, let $\mu_0(x)$ be a given L-periodic function in
$L^\infty(\R^N)$ taking only the two values $\mu_0^+=B$ and
$\mu_0^-=0$. Let $\lambda_{1,0}$ be the first eigenvalue of the
operator $-\nabla^2-\mu_0 I$ on $C$, with L-periodicity conditions,
$\phi_0$ the associated eigenfunction with minimal value
$\underline{\phi_0}$, and
\[
\delta_{1,0}:=\ds{\frac{\lambda_{1,0}^2
\underline{\phi_0}}{(1+\underline{\phi_0})^2}} \quad\hbox{and}\quad
\delta_{2,0}:=\ds{\frac{\lambda_{1,0}^2 }{4}}.
\]
We have the following proposition.

\begin{a3proposition}\label{propgap}
Assume that $\mu(x)=\mu_0(x)+\mu^-,$ with $\mu^->\lambda_{1,0}$. Let
$\delta_1$ and $\delta_2$ be defined by {\rm(\ref{defd1d2})}. Then
we have $\delta_2-\delta_1=(1-\frac{\mu^-}{\lambda_{1,0}})^2
(\delta_{2,0}-\delta_{1,0})$.
\end{a3proposition}

This result also indicates that the information on $\delta^{\ast}$
is all the more precise as the growth rate function takes low
values. However, the ``relative thickness'' of the interval
$(\delta_1,\delta_2)$, compared to $\delta_1,
\frac{\delta_2-\delta_1}{\delta_1}$, does not depend on $\mu^-$, as
can be easily seen.

{\it Proof of Proposition\/ {\rm\ref{propgap}}}. The relation
$\lambda_1[\mu(x)]=\lambda_{1,0}-\mu^-$ is a direct consequence of
the uniqueness of the first eigenvalue $\lambda_1$. We assume that
$\mu^->\lambda_{1,0}$, so that $\lambda_1[\mu(x)]<0$. From the
uniqueness of the eigenfunction $\phi$ associated with $\lambda_1$,
$\phi$ does not depend on $\mu^-$. Therefore, $\delta_1$ and
$\delta_2$ satisfy $\delta_1= \frac{(\lambda_{1,0}-\mu^-)^2
\underline{\phi_0}}{(1+\underline{\phi_0})^2}$ and $\delta_2=
\frac{(\lambda_{1,0}-\mu^-)^2 }{4}$. The result immediately follows.
\qquad $\Box$

\section{A few comments on the proportional harvesting model}\label{sec_prop}

In this model, the population density $u$ is governed by the
equation \be\label{eqevoprop} u_t=D \nabla^2 u+u(\mu(x)-\nu(x)
u)-q(x) u, \quad x\in\Omega, \ee with L-periodicity of the functions
$\mu(x)$, $\nu(x)$, and $q(x)$ in the periodic case, and with
Neumann or Dirichlet boundary conditions in the bounded case.
Setting
\[
\tau(x):=\mu(x)-q(x),
\]
this model becomes equivalent to the SKT model (\ref{eqi4}). Hence,
many properties of the solutions of this model are described in the
existing literature. In particular the existence, nonexistence, and
uniqueness results of Theorems \ref{thbhr1} and~\ref{thdbhr1} apply.
The condition $\lambda_1[\mu(x)-q(x)]<0$ is therefore  necessary and
sufficient for species persistence. Furthermore, the theoretical
results of \cite{bhr1}, \cite{ccL},  \cite{rh1}, \cite{rs1} on the
effects of habitat arrangement on species persistence are also true
for this model.

For instance, when the function $\mu(x)$ is constant, with
$\mu(x)\equiv\mu_1>0$, and if the domain $\Omega$ is convex and
symmetric with respect to each axis $\{x_1=0\},\ldots,\{x_N=0\}$,
the next result is a straightforward consequence of the
paper~\cite{bhr1}.

\begin{a3theorem}
{\rm(i)} In the periodic case{\em,} $\lambda_1[\mu_1-q_k^*(x)]\leq
\lambda_1[\mu_1-q(x)]$.

{\rm(ii)} In the bounded Dirichlet case{\em,}
$\lambda_1[\mu_1-q_k^*(x)]\leq \lambda_1[\mu_1-q(x)]$.

{\rm(iii)} In the bounded Neumann case{\em,} if\/ $\Omega$ is a
rectangle, $\lambda_1[\mu_1-q_k^{\sharp}(x)]\leq \lambda_1[\mu_1-
q(x)]$.\label{th3}
\end{a3theorem}

Here $q_k^*$ denotes the symmetric decreasing Steiner rearrangement
of the function $q$ with respect to the variable $x_k$, and
$q_k^{\sharp}$ denotes the monotone rearrangement of $q$ with
respect to $x_k$ (see \cite{bhr1} and~\cite{blr} for the definition
of these rearrangements). These rearrangements of a function $q$
preserve not only its mean value, but also its distribution
function. This means that if, for instance, $q$ corresponds to a
``patch'' function taking the values $q_{1}$, $q_{2}$, and $q_{3}$
in some regions $A_{1}$, $A_{2}$, and $A_{3}$, respectively, with
$A_{1}+A_{2}+A_{3}=|C|$, then the areas of the regions where the
rearranged functions $q^*$ and $q^\sharp$ take the values $q_{1}$,
$q_{2}$, and $q_{3}$ remain equal to  $A_{1}$, $A_{2}$, and $A_{3}$,
respectively.

\begin{figure}[t!]
\centering
\subfigure[]{\label{figa1}\includegraphics[width=5cm]{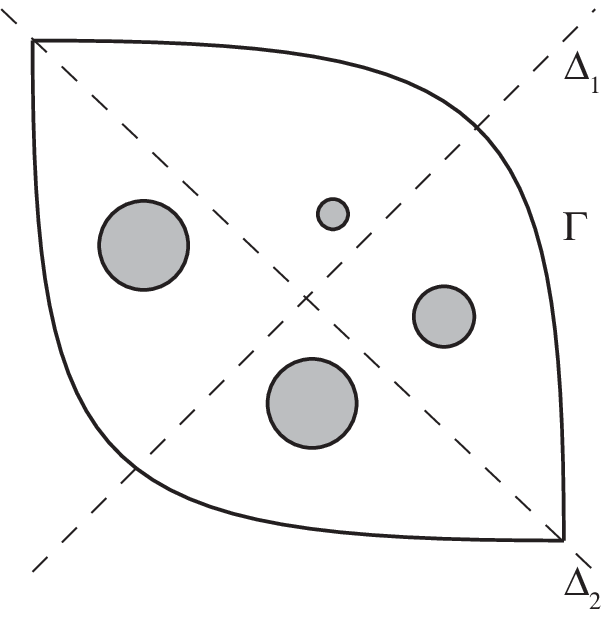}}
\hspace{1cm}
\subfigure[]{\label{figa2}\includegraphics[width=5cm]{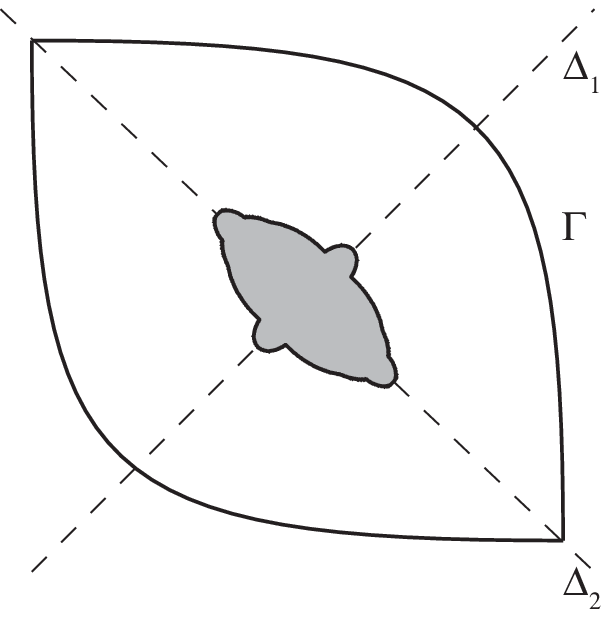}}
\subfigure[]{\label{figa3}\includegraphics[width=5cm]{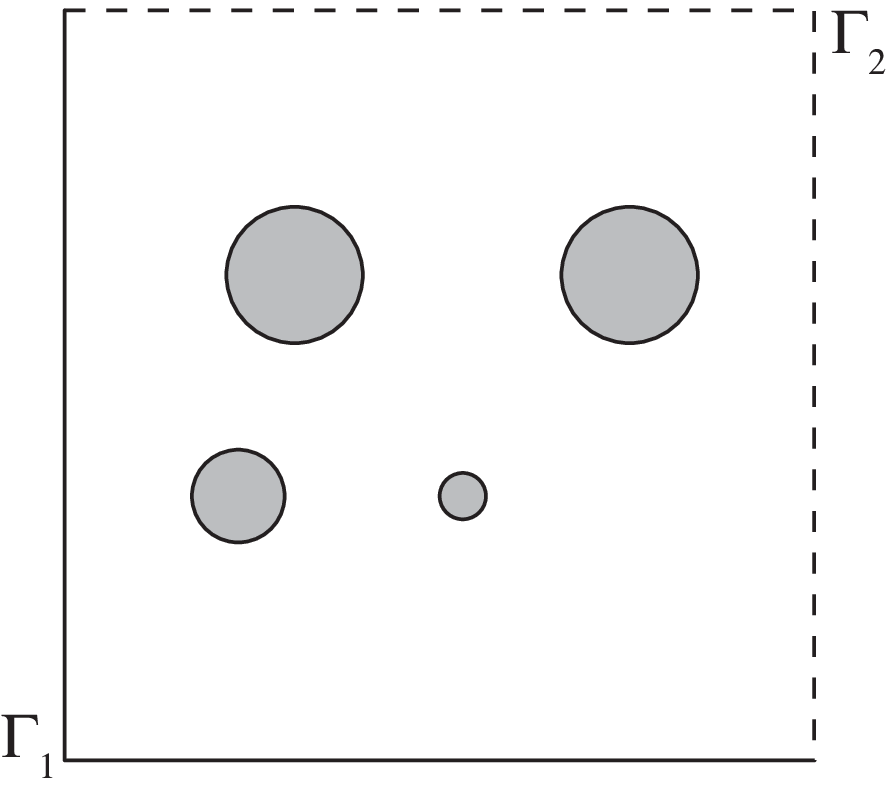}}
\hspace{1cm}
\subfigure[]{\label{figa4}\includegraphics[width=5cm]{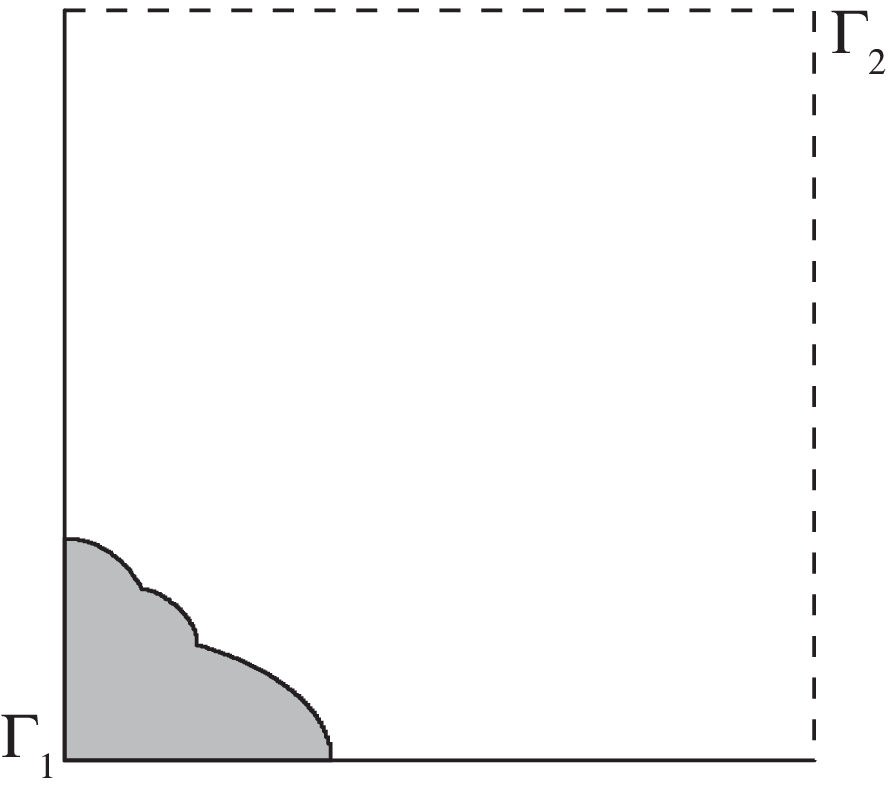}}\vspace*{-6pt}
\caption{Examples of applications of Theorem~{\rm\ref{th3}.ii--iii}
to reserves management. In panels {\rm(a)} and~{\rm(b),} the
boundary $\Gamma$ of $\Omega$ is lethal (Dirichlet boundary
conditions). {\rm(a)} The initial effort function $q(x)$  takes two
values{\em,} $q^+>0$ in the white area, and $q^-=0$ in the shadowed
regions{\em,} which correspond to reserves. {\rm(b)} Position of the
reserves after a symmetric decreasing Steiner rearrangement along
the $\Delta_1$ and $\Delta_2$ axes{\em,}
 successively. The rearranged configuration~{\rm(b)} always give more chances
of species persistence.  In panels {\rm(c)} and~{\rm(d),} the
boundary $\Gamma$ is divided into two parts{\em:} $\Gamma=\Gamma_1
\cup \Gamma_2$. $\Gamma_1$ is represented with a solid line and can
correspond to a coast{\em,} while $\Gamma_2$ is represented with a
dashed line and can correspond to a nonphysical limit that the
species cannot cross (Neumann boundary conditions). {\rm(c)} The
effort function $q(x)$ again takes two values{\em,} $q^+>0$ in the
white area, and $q^-=0$ in the reserves. {\rm(d)} Position of the
reserves after monotone rearrangement along  the horizontal and
vertical axes{\em,} successively. The chances of persistence are
better in the rearranged
configuration~{\rm(d)}.}\label{fig:exprop}\vspace{-1.5pc}
\end{figure}

Theorem~\ref{th3} combined with Theorem~\ref{thdbhr1} says that the
spatially rearranged harvesting strategies are better for species
survival. This result can be helpful from a resource management
point of view. Indeed, the authorities can rearrange the position of
the harvested areas in order to improve the chances of population
persistence. The result of Theorem~\ref{th3} shows that, in the
framework of these models, the creation of a large reserve gives
persistence more chances than the creation of several small
reserves, and is in accordance with the former results of
\cite{shi1}~and~\cite{neub} in the Dirichlet case. See
Figure~\ref{fig:exprop} for some illustrations in the bounded case
with Dirichlet and Neumann boundary conditions.

\section{Discussion}
We have proposed a model for the study of populations in
heterogeneous environments, for populations submitted to an external
negative forcing term. This forcing term could be regarded as a
``quasi-constant-yield'' harvesting, \hbox{depending} only on the
population density $u$ when $u$ is below a certain small threshold
$\varepsilon$. The introduction of such a threshold $\varepsilon$
was necessary for ensuring the nonnegativity of the solutions of our
model, and therefore its actuality.

We carried out new mathematical results on the elliptic equation
satisfied by the stationary states of the model, and on the
associated parabolic equation. Both qualitative and quantitative
results were obtained.

From the qualitative point of view, we described the behavior of the
model solutions in terms of the harvesting amplitude $\delta$. Two
main types of stationary solutions were found: the remnant
solutions, always below a small threshold $\varepsilon_0$ and
therefore close to~0, and the significant solutions,  always above
this threshold, thus ensuring a time-constant yield. We discussed
the maximum number of significant stationary solutions, which we
found equal to~2, under a hypothesis of positivity of the second
eigenvalue $\lambda_2$ of a linear operator. We further investigated
the long-time behavior of the solution of our model, starting from a
nonharvested population at equilibrium. We found a critical value
$\delta^*$ of the harvesting term amplitude, below which the
population density tends over time to a significant stationary
solution, and above which it converges to a stationary solution
which is not significant. We also established quantitative formulae
for some lower and upper bounds for $\delta^*$: $\delta_{1}$ and
$\delta_{2}$, respectively. The threshold $\delta_2$ has the
additional property that, whenever the amplitude $\delta$ is above
$\delta_{2}$, the population density decreases to a remnant
stationary solution.

The quantitative aspects of our study mainly consisted of discussing
the effect of environmental fragmentation on these thresholds
$\delta_1$ and $\delta_2$, and therefore on the interactions between
environmental fragmentation and maximum sustainable yield. Namely,
when computing the values of $\delta_{1}$ and $\delta_{2}$ on 2000
samples of stochastically obtained patchy environments, with
different levels of fragmentation, we found an \hbox{increasing}
relationship between these two coefficients and an environmental
aggregation index $s$. This indicates that, for given areas of
favorable and unfavorable regions, the harvesting quota that a
species can sustain, while ensuring a time-constant yield, is higher
when the favorable regions are aggregated.

The reader may note that, in our model,  the species mobility was
not affected by the environmental heterogeneity. Such a dependence
could be modeled by using  a more general dispersion term, of the
form $\nabla\cdot (A(x)\nabla u)$, instead of $D \nabla^{2}u$, where
$A(x)$ stands for the diffusion matrix (see \cite{bhr1}, \cite{sk}).
In fact, most of our results still work when the matrix $A$ is of
class $C^{1,\alpha}$ (with $\alpha>0$) and uniformly elliptic,
i.e.,~when there exists $\tau>0$ such that $A(x)\ge\tau I_N$ for all
$x\in\Omega$. Indeed, Theorems \ref{thbhr1}, \ref{thdbhr1},
\ref{th1.1}, \ref{th1}, and \ref{th2} remain true under this more
general assumption. However, the effects of environmental
heterogeneity may differ, depending on the way $A(x)$ and $\mu(x)$
are correlated (see~\cite{kks}). In the proportional harvesting
case, the results of section~4 on the effects of the arrangements of
the harvested regions may also not be valid with this dispersion
term. However, in situations where $A(x)$ takes low values (slow
motion) when $q(x)$  is low (``reserves''; see section~4), as
underlined in~\cite{rh1}, a simultaneous rearrangement of the
functions $A(x)$ and $q(x)$ would lead to lower $\lambda_{1}$ values
and therefore to higher chances of species survival.\vspace*{-2pt}

\section*{Acknowledgment}

The authors would like to thank the anonymous referees for their
valuable suggestions and insightful comments.

\end{document}